\documentclass[twoside]{amsart}
\newtheorem{thm}{Th\'eor\`eme}[section]

\newtheorem{prop}[thm]{Proposition}

\newtheorem{lem}[thm]{Lemme}
\newtheorem{dfn}[thm]{D\'efinition}
\newtheorem{remark}{Remarque} 

\begin{document}
{
\center{\Large{\bf 
Connexion de Gau\ss - Manin associ\'ee \`a la d\'eformation 
verselle de la singularit\'e $A_\mu$
et z\'eros de l'int\'egrale hyperelliptique.
} }
 \vspace{1pc}
\center{\large{ Susumu TANAB\'E $^{1}$}}
 }
 \vspace{0.5pc}

\begin{center}
\begin{minipage}[t]{12.2cm}
{\sc R\'esum\'e -} {\em
On \'etudie le syst\`eme des \'equations differ\'entielles satisfaites
par l'int\'egrale hyperelliptique  associ\'ee au cycle
 $\gamma_s \subset \{(x,y) \in {\bf R}^2:
H(x,y;s)
=0 \}$ d\'efinie pour la d\'eformation verselle de la singularit\'e $A_\mu.$
Comme application,
on obtient une estimation de la multiplicit\'e  des z\'eros  de
l'int\'egrale $I_{\omega}(s)=
\int_{\gamma_s} \omega $ en fonction de $\mu$ et de $deg (\omega).$
 }
\end{minipage}
\end{center}

{
\center{\section
{\bf Introduction }}
}

Dans cette note on poursuit directement les recherches sur les 
probl\`emes trait\'es 
dans la premi\`ere partie du travail pr\'ecedent \cite{AT2}, i.e. la 
description raisonnable de la connexion de Gauss-Manin associ\'ee \`a la d\'eformation
verselle de la singularit\'e $A_\mu$.
Notre objet principal est 
l'int\'egrale hyperelliptique que l'on 
d\'efinit sur une courbe hyperelliptique. 
Regardons un polyn\^ome d\'ependant de $\mu -1$ param\`etres
$ s'= (s_1,$ $s_2,$ $\cdots, s_{\mu -1}),$
$$ H(x,y) = F(x,s') - y^{2}$$
$$ F(x,s')=x^{\mu +1} + s_{\mu -1}
x^{\mu -1} + \ldots + s_1 x ,\;\;\;\; \mu, \nu \geq 2.  \leqno(1.1) $$
On prend un cycle \' evanescent 
$$\gamma_s \subset
\{(x,y)\in {\bf R}^2: H(x,y)+ s_0=0 \}, \;s=(s_0,s') \leqno(1.2)$$ du  polyn\^ome $ H(x,y) $ et une 
1-forme polynomiale $\omega= P (x,y) dx + Q(x,y) dy,$ alors
on
appelle  l'expression  $$  I_{\omega}   (s)   =   \int_{\gamma_s}
\omega,$$ int\'egrale hyperelliptique.

Dans \cite{AT2}, on s'occupait, entre autres,  du probl\`eme 
d'adjacence entre diverses 
int\'egrales hyperelliptiques associ\'ees aux differentes
singularit\'es $A_{\mu-1}$ et $A_{\mu}$. Par contre, ici on recherche
la structure du syst\`eme de Gauss-Manin le long 
d'une strate de son ensemble critique 
(du discriminant). L'invariance des exposants de monodromie le long de la 
strate
$\mu = const.$ est connue depuis Varchenko \cite{Var2}. Ici nous allons 
\'etablir un \'enonc\'e un peu plus fin sur les exposants 
caract\'eristiques du syst\`eme de Fuchs (voir la d\'efinition ~\ref{dfn31}) qui d\'ependent explicitement du 
degr\'e de l'int\'egrand $\omega$ de $  I_{\omega}(s)$ (voir le th\'eor\`eme
d'isomonodromie renforc\'e ~\ref{thm3}). Pour autant qu'il s'agisse de 
l'int\'egrale hyperelliptique g\'en\'eralis\'ee, notre th\'eor\`eme est une
version forte de celui de Varchenko, car nous constatons l'invariance de 
$\mu$ (ou bien $\mu +1$ ) exposants pour chaque int\'egrale 
$  I_{\omega}   (s),$ pourtant Varchenko a d\'emontr\'e l'invariance
du minimum de ces $\mu$ exposants. 

Dans \cite{AT1}, \cite{AT2} (17), nous avons propos\'e une nouvelle
d\'efinition du syst\`eme du type de Fuchs avec lieu singulier $D,$
un diviseur d'une vari\'et\'e complexe lisse $S$
en tant que syst\`eme de Pfaff avec les coefficients
de $\Omega_S^1 (log\; D),$ formes diff\'erentielles logarithmiques. Il est 
naturel de se poser la question de savoir comment le syst\`eme du type de Fuchs
ainsi d\'efini se lie \`a l'\'equation de Fuchs au sens classique? 
Le th\'eor\`eme ~\ref{thm3} r\'epond, entre autres, \`a cette question aussi.Je
tiens \`a noter qu'une formule hypoth\'etique propos\'ee par V.P.Palamodov
servait de probl\`eme moteur de notre recherche. Il se demande quel 
op\'erateur diff\'erentiel doit annuler  $I_{dx}(t).$ Notre proposition
 2.4 fournit une reponse.

 \vspace{1.5pc}
\footnoterule

\footnotesize{AMS Subject Classification: primaire 34C08, 14K20, 14D05,
secondaire 32S30, 33C20, 34M99.

Mots cl\'es: connexion de Gauss-Manin, cycles limites,
l'\'equation du type de Fuchs, d\'eformation isomonodromique.}

 \footnotesize{${}^1$  Travail   r\'ealis\'e   par   le   soutien
financier de l'
homme d'affaires M.Mikhail S.Gavounas (Moscou, Russie)
et du Max Planck Institut f\"ur Mathematik, Bonn}

\normalsize
\newpage
\par

A partir du chapitre 4, notre pr\'eoccupation sera l'\'etablissement 
d'une estimation par haut de la multiplicit\'e des z\'eros de 
l'int\'egrale hyperelliptique. C'est une demarche vers une
r\'eponse raisonnable au XVIe probl\`eme de Hilbert sur le nombre des cycles
limites.
On prend un 
hamiltonien polynomial comme (1.1) (ou bien plus g\'n\'eralement $ H(x,y) = F(x,s') - y^{\nu}$), et impose la condition que 
pour chaque valeur critique $s_0^{(i)}
(s_0^{(i)}\not = s_0^{(j)}, $ si $i \not = j$)
$ F(x,s')+s_0^{(i)}=0$
ait singularit\'e du type de Morse i.e. le hessien soit non-d\'eg\'en\'er\'e
\`a chaque point critique.

Regardons       un       polyn\^ome
$$P_{K,m}(x,y)        =
\prod_{i=1}^{L}(x-x^{(i)})^{k_i} y^m  , \leqno(1.3)$$
$ k_i \in {\bf N}, m \in {\bf Z}$ tel que $ x^{(2)}, \cdots, x^{(L)} $
ne sont pas de  points  critiques  de  $F(x,s').  $  Par  contre,
 $x^{(1)}$ peut \^etre un point critique de
$F(x,s').$ On note $K= \sum_{i=1}^L k_i.$

Puisque nous visons  \`a \'etablir l'estimation de la multiplicit\'e des 
z\'eros 
de $  I_{\omega}(s),$ pr\'ecisons la d\'efinition de cette notion.

\begin{dfn} {\em
Si une s\'erie convergente pr\`es de $t=t_0$ d\'efinit une fonction multivalu\'ee dans un secteur $\Delta_\theta = \{t : |t-t_0| < \theta\} \subset {\bf C},$
$$ f(t) = \sum_{\rho \in {\bf Q}, k \in {\bf N} }f_{\rho,k}(t-t_0)^\rho (log(t-t_0))^k$$
cette s\'erie s'appelle fonction de Dulac.

Nous disons que une fonction de Dulac a $t_0$ comme z\'ero de 
multiplicit\'e $(k_1+1)([\rho_1] +1)$  si 
$$  \rho_1 = \{ min\; \rho; f_{\rho,k} \not =0 \}, $$
$$  k_1 = \{ max \; k; f_{\rho_1,k} \not =0 \}. $$}
\label{dfn11}
\end{dfn}
D'ici bas, la notation $[\rho]$ signifie la partie enti\`ere d'un nombre 
rationel $\rho.$
Tout au long de cet article, on comprendra la multiplicit\'e des
z\'eros de $I_{\omega}(s)$ comme celle des z\'eros d'une fonction de Dulac en 
variable $s_0.$

\begin{thm} Dans la situation ci-dessus, on consid\`ere
l'int\'egrale
hyperelliptique 
$I_{P_{K, m}}(s)$
prise le long d'un cycle $\gamma_s =
\{(x,y)\in {\bf R}^2: H(x,y)+ s_0=0 \},$
$$   I_{P_{K,m}}(s)=   \int_{\gamma_s}P_{K,m}(x,y)    dx    ,
\leqno(1.4)$$
avec $ K \in {\bf N}, m \in {\bf Z}.$ Supposons $ I_{P_{K, m}}(s)\not \equiv
0.$  Alors on a r\'esultats suivants.

i) Si $\mu$ pair, la multiplicit\'e $N$ des z\'eros
de  
l'int\'egrale $I_{P_{K, m}}(s)$ \`a l'un des points de ramification
${\tilde t}\in \{ s_0^{(1)}, \cdots ,s_0^{(\mu)}    \}$ v\'erifie:
$$ N \leq 
2[\frac{K+m+\mu}{2}]. \leqno(1.5)$$

ii) Supposons que $k_1=0$ dans l'expression (1.3).
Alors la multiplicit\'e $N$ des z\'eros
de  l'int\'egrale $I_{P_{K, m}}(s)$ avec $\mu$ impair, 
\`a l'un des points de ramification
${\tilde s_0}\in \{ s_0^{(1)}, \cdots ,s_0^{(\mu)}    \}$ 
v\'erifie:
$$ N \leq max \{ \mu-1, 
2[\frac{m+3}{2}] \}.
\leqno(1.6)$$

iii) La multiplicit\'e $N$ des z\'eros
de  
l'int\'egrale $I_{P_{K, m}}(t)$ au point 
$\tilde s_0 \not \in \{ s_0^{(1)}, \cdots ,s_0^{(\mu)}\}$ ne 
d\'epasse pas $\mu +K .$
\label{thm0}
\end{thm}

Pour \'etablir le th\'eor\`eme, dans un contexte plus g\'en\'eral que 
celui de l'int\'egrale hyperelliptique, nous allons \'etudier l'int\'egrale 
d'Abel associ\'ee \`a une courbe d\'efinie par un polyn\^ome 
$$ H(x,y)=x^{\mu +1} + s_{\mu -1}
x^{\mu -1} + \ldots + s_1 x -y^\nu ,\;\;\;\; \mu, \nu \geq 2.  \leqno(1.1)' $$
Si on prend un cycle \' evanescent 
$$\gamma_s \subset
\{(x,y)\in {\bf R}^2: H(x,y)+ s_0=0 \}, \;s=(s_0,s') \leqno(1.2)'$$ 
de la courbe   et une 
1-forme polynomiale $\omega= P (x,y) dx + Q(x,y) dy,$ alors il est 
possible de d'efinir l'int'egrale d'Abel
$$  I_{\omega}   (s)   =   \int_{\gamma_s}
\omega,$$ d'une fa\c{c}on analogue \`a l'int\'egrale hyperelliptique. 
Voir Th\'eor\`eme {thm51}.

Il est facile de d\'eduire de $(1.4), (1.5)$ que si 
$\mu, K ,m \leq n, $
alors la multiplicit\'e des z\'eros $N$ est domin\'e 
par une fonction lin\'eaire 
de $n.$ 
Nos r\'esultats, donc,  paraissent une  am\'elioration  des
estimations obtenues dans   \cite{Mar} ( la multiplicit\'e de chaque
z\'ero $ \leq \frac{n^4 + n^2 -2}{2}$  )  
pour autant qu'il s'agisse de l'int\'egrale du type  $(1.4).$  

L'auteur tient \`a remercier L. Gavrilov qui m'a donn\'e des cours
d'initiation  \`a ce sujet,
S. Yakovenko, J. P. Fran\c{c}oise, P.Mardesic et D. Novikov de 
leurs
conseils gentils.

 \vspace{2pc}

{\center{\section{\bf Les \'enonc\'es sur le syst\`eme de Gauss-Manin  associ\'e
aux singularit\'es $A_\mu$ }}}

 On se souvient ici des r\'esultats principaux de \cite{AT2} concernant
les int\'egrales d'Abel 
$$K^{\lambda}_i(s) = \int z^i (z^{\mu +1} + s_{\mu -1} z^{\mu -1}+ \ldots
+ s_1 z +s_0 )^\lambda dz, \hspace{2pc} i = 0, \ldots, \mu +1. \leqno(2.1) $$
L'int\'egrale $(2.1)$ est un prototype  de notre recherche,
pour
un hamiltonien du type $(1.1)'. $
Notre  d\'emarche  consiste   en   l'analyse   de l'int\'egrale
$(2.1)$ \`a l'aide de l'op\'erateur  diff\'erentiel
qui l'annule.

Regardons   les   int\'egrales  de  p\'eriodes   pour   une
vari\'et\'e alg\'ebrique de dimension complexe un
$ X^{(\mu,\nu)}_s= \{(z,y)\in P^2{\bf C}; y^{\nu} - F(z,s')= s_0\}$
param\'etr\'ee par $\mu$ param\`etres $(s_0,$ $ \ldots,
s_{\mu-1})$ avec
$$  F(z,s') = z^{\mu + 1} +s_{\mu-1} z^{\mu -1}+ \ldots + s_1 z.$$
Le polyn\^ome $F(z,s')+s_0$ donne une d\'eformation verselle de la singularit\'e $A_\mu,$ $  F(z,0) = z^{\mu + 1}.$
Remarquons que $ rang H^1(X^{(\mu,\nu)}_s)=\mu(\nu -1),  $
et comme base de $H^1(X^{(\mu,\nu)}_s),$ on peut choisir
$ x^k y^{\ell}dx, ( 0 \leq k \leq \mu -1, 1 \leq \ell \leq \nu -1  ).$
Il est facile de voir que l'action de quasihomog\'en\'eit\'e agit sur 
$  F(z,s') +s_0:$
$$ \tau : (s_0, \cdots, s_{\mu -1},z) \rightarrow
( t^{\mu  +1}  s_0,  t^{\mu   }
s_1,\cdots,   t^2   s_{\mu -1}, tz),  \;\; t \in {\bf R}_+.$$

Etant fix\'e un cycle \'evanescent $ \gamma_s \in H_1(X^{(\mu,\nu)}_s),$ 
on  consid\`ere l'int\'egrale  de p\'eriodes
comme celle d\'efinie  le long d'un cycle  $
Reg(\gamma_s)$  sur  $\tilde{\bf
C}_x$ un rev\^etement de $\bf C$ avec $\nu$ feuilles.
Si on note  $\lambda = \frac{\ell}{\nu}:$
$$I_{x^iy^\ell dx,\gamma_s}(s)=K^{\lambda}_{i, \gamma_s}(s) = 
\int_{ Reg(\gamma_s)} z^i (F(z,s') + s_0)^\lambda dz,
\hspace{2pc}
i = 0,  \ldots,\mu  +  1. \leqno(2.2) $$
Ici $ Reg(\gamma_s)$ d\'enote un cycle appel\'e r\'egularis\'e, de
 sorte que
l'int\'egration le long de celui-ci soit bien d\'efinie. Sur tout il doit
\^etre choisi de fa\c{c}on que $K^{\frac{\ell}{\nu}}_{i, \gamma_s} (s)$ 
soit aux
valeurs r\'eelles pour $s\in {\bf R}$ et $\gamma_s$ un cycle \'evanescent
r\'eel. Il est obtenu en
"gonflant" $\gamma_s$ \`a l'aide de l'op\'erateur de Leray (voir 4.2.
\cite{Vas}).
Dans notre situation $  Reg(\gamma_s)$  n'est  qu'une  somme  des
lacets doubles
de Pochhammer \cite{AK}. Gr\^ace au th\'eor\`eme des r\'esidus de Leray
\cite{Br}, on n'a pas besoin de se soucier de concr\'etiser le cycle
$ Reg(\gamma_s),$ lors de l'\'etablissement des \'equations diff\'erentielles
satisfaites par l'int\'egrale $ (2.2).$ Toutes les int\'egrales d\'efinies par
$(2.2)$ satisfont la m\^eme \'equation diff\'erentielle ind\'ependante
du cycle $\gamma_s.$ Donc on \'ecrira souvent $K^{\lambda}_{i}(s)$
au lieu de $K^{\lambda}_{i, \gamma_s}(s)$ sinon on a besoin de pr\'eciser
le cycle d'int\'egration.
En    qu\^ete   d'une    expression    concr\`ete    des    int\'egrales
de p\'eriodes autour de ses  points de ramification,
nous partons de la proposition suivante.

\begin{prop}[\cite{M}, \cite{P}]
Les int\'egrales de p\'eriodes $K^{\lambda}_0(s),$
$\ldots,$
$K^{\lambda}_{\mu + 1}(s)$ de (2.2)  satisfont    le    syst\`eme
holon\^ome suivant d'\'equations diff\'erentielles:
$$\sum_{\ell =0}^{\mu -1 }
s_{\ell} \frac{\partial}{\partial s_0}
K^{\lambda}_{\ell+i} + \frac {\partial}{\partial s_0}
K^{\lambda}_{\mu +1+i} = {\lambda}K^{\lambda}_i  ,
\;\;\; 0 \leq i \leq \mu-1,\leqno(2.3)_i$$
$$\sum_{\ell =1}^{\mu -1} {\ell}s_{\ell} \frac{\partial}{\partial s_0}
K^{\lambda}_{\ell+j} +(\mu +1) \frac{\partial}{\partial s_0}
K^\lambda_{\mu +1+j} = -(j+1)K^{\lambda}_j  ,
\;\;\; -1\leq j \leq \mu -1. \leqno(2.4)_j$$
Nous avons une repr\'esentation matricielle entre les
int\'egrales:
$$ {\bf \Sigma} \cdot \vec{b} = \vec{a},$$
o\`u
\begin{center}
$ {\bf \Sigma}=
\tiny{ \left [
\begin{array}{ccccccccccc}
s_0&s_1&\cdots&s_{\mu -2} & s_{\mu-1}&0&1&\cdots&0&0&0\\
0&s_0&\cdots &s_{\mu-3}&s_{\mu -2} & s_{\mu-1}&0&\cdots&0&0&0\\
0&0&\cdots&s_{\mu -4} &s_{\mu-3}&s_{\mu -2} & s_{\mu-1}&\ddots&0&0&0\\
\vdots&\vdots &\ddots&\vdots&\vdots &\vdots&\vdots &\cdots& \vdots&\vdots &\vdots \\
0&0&\cdots&s_1&s_2 &s_3&s_4 &\cdots& 1 &0 & 0 \\
0&0&\cdots&s_0&s_1&s_2&s_3 &\cdots& 0&1 & 0 \\
0&0&\cdots&0&s_0&s_1&s_2& \cdots &s_{\mu}&0&1 \\
s_1&2s_2&\cdots&(\mu -1)s_{\mu-1} &0&\mu +1&0&\cdots&0&0&0\\
0&s_1&\cdots&(\mu-2)s_{\mu-2}
&({\mu} -1)s_{\mu-1} &0&{\mu +1}&\ddots&0&0&0\\
\vdots&\vdots &\cdots&\vdots&\vdots &\vdots&\ddots &\cdots&\vdots
&\vdots&\vdots\\
0&0&\cdots&s_1&2s_2 &3s_3&\cdots &\cdots& \mu +1 &0 & 0 \\
0&0 &\cdots
&0&s_1&2s_2&\cdots&\cdots&0&{\mu +1}&0\\
0&0 &\cdots&0
&0&s_1&\cdots&\cdots&({\mu}-1)s_{\mu-1}&0&{\mu +1}\\
\end{array}
\right] } ,$
\end{center}
$$\vec{a}= ^t({\lambda}K^{\lambda} _{0},{\lambda}K^{\lambda}_{1 },\cdots,
{\lambda}K^{\lambda}_{\mu-1},0,
-K^{\lambda}_{0}, -2K^{\lambda}_{1},\cdots, -\mu K^{\lambda}_{\mu -1} )  $$
$$ \vec{b}= ^t (\frac{\partial}{\partial s_0}K^{\lambda} _{0},
\frac{\partial}{\partial s_0}K^{\lambda}_{1},\cdots,
\frac{\partial}{\partial s_0}K^{\lambda}_{2\mu }).$$
Remarque:
$\frac{\partial}{\partial s_i}K^{\lambda}_{j}(s)=
\frac{\partial}{\partial s_j}K^{\lambda}_{i} (s).$

\label{prop1}

\end{prop}

Effectivement, $2\mu$ int\'egrales de p\'eriodes prennent part aux
\'equations $(2.3), (2.4),$ au lieu  de  $\mu$  int\'egrales.  En
les supprimant, on obtient les relations  syzygy  non-triviales
entre $\mu$ int\'egrales de p\'eriodes.

\begin{prop}(\cite{AT2})
Les int\'egrales  ${\bf K}(s) =$ $^t( K_0(s), \cdots, $ $ K_{\mu -1}(s)),$
$s = ( s_0, \cdots,$ $ s_{\mu -1} )$ satisfont le  syst\`eme  holon\^ome
suivant d'\'equations diff\'erentielles:
$$
S\frac {\partial}{\partial s_0}{\bf K} =(L+V(s_2,\cdots,s_{\mu -1})){\bf K},
\leqno(2.5)$$
ou $ V(s') = $
\begin{center}
$  = {\displaystyle \frac{1}{(\mu +1)^2}}
\left[
\begin{array}{ccccccc}
0&0 & 0 & 0 & \cdots& 0 & 0 \\
0&0 & 0 & 0 & \cdots &0 & 0 \\
2s_{\mu-1}&0 & 0 & 0 &\cdots  & 0 & 0 \\
3s_{\mu -2}& 2 \cdot 2 s_{\mu -1} & 0 & 0 & \cdots & \vdots & \vdots \\
4s_{\mu -3}& 3 \cdot 2 s_{\mu -2} & 2 \cdot 3 s_{\mu -1} & 0 &\cdots  &
\vdots & \vdots\\
\vdots &\vdots & \vdots & \ddots & \cdots & \vdots & \vdots \\
(\mu -1)s_2& 2(\mu -2 )s_3 & 3(\mu -3)s_4 & \cdots & (\mu -2)2s_{\mu -1}
& 0 & 0 \\
\end{array}
\right],$
\end{center}
et $ s'=( s_1 ,$ $ \cdots ,$ $ s_{\mu -1} ) . $
Les \'elements $v_{i,j}$ de la matrice $V(s_2,\cdots,s_{\mu -1}) $
sont d\'etermin\'es par la r\`egle ci-dessous,
$$(j+1) v_{i,j} = j v_{i+1,j+1}, \,  v_{1,j}=(j-1)s_{\mu -j + 2}, \;\;
1 \leq i \leq \mu -2, \, 3 \leq j \leq \mu.$$
La matrice $S$ admet une \'ecriture comme suit,
$$S = s_0 {\rm id}_{\mu} + C(s'), $$
avec une matrice polynomiale  $C(s')$ et une
matrice diagonale $L  $ repr\'esentant les poids quasihomog\`enes de
formes diff\'erentielles correspondants \`a $K_j^{\lambda}(s)$,
$$ L = {\rm diag} (\lambda + \frac{1}{\mu +1} ,\ldots, \lambda +
\frac{\mu}{\mu +1}).$$
\label{prop2}
\end{prop}

Nous nous servons de  la notation  $\Delta_{\mu}(s)$  d\'esignant 
un  polyn\^ome monique en la variable $s_0,$ qui peut \^etre consid\'er\'e
comme  le
 discriminant du polyn\^ome $z^{\mu+ 1} + s_{\mu-1} z^{\mu -1} + \ldots
+ s_1 z +s_0. $  Il est calcul\'e par la matrice $S(s)$ de la
Proposition ~\ref{prop2},
$$ \Delta_{\mu}(s_0,s') = \det \;S(s) = \det \Sigma .$$

\begin{dfn}
On dit que $s' \in {\bf R}^{\mu -1}$ appartient \`a l'ensemble de bifurcation
$B \subset {\bf R}^{\mu -1}$ si et seulement s'il existe $s_0 \in {\bf C}$
tel que $$\Delta_{\mu}(s_0,s') = \frac{\partial}{\partial s_0}\Delta_{\mu}(s_0,s')=0.$$
\label{dfn}
\end{dfn}

Si  $s' \in B$ l'\'equation $F(z,s') =t$ poss\`ede soit des racines mutiples
soit des valeurs critiques multiples 
pour une certaine valeur de $''t.''$ L'ensemble $B$ est une vari\'et\'e
alg\'ebrique dans ${\bf R}^{\mu -1}$ de codimension 1.

 Nous nous rappelons ici une relation  de la matrice  $S$  avec  le
champ de vecteurs logarithmiques formul\'e par K.Saito
(\cite{S3}).

\begin{lem}
Si  on  utilise la notation
$ S(s)$ $ = $ $(\sigma_{i,j}(s))_{0\leq i,j \leq \mu -1},$ alors les
vecteurs $\xi_i, \; i=0,\cdots, \mu -1$ d\'efinis comme suit
$$    \xi_i    =    \sum_{j=0}^{\mu     -1}     \sigma_{i,j}\frac
{\partial}{\partial s_j} $$
constituent  le  champ  de  vecteurs  logarithmiques  tangent  au
discriminant  $D  =  \{\Delta_{\mu}(s)  =  \det   \;S(s)   =0\}.$
Autrement dit, $\xi_0,\cdots, \xi_{\mu-1}$  forment  une  base
libre de $Der_{{\bf  C}^{\mu}}(\log  D)$  en  tant  que
${\mathcal O}_{{\bf  C}^{\mu}}- $ module.
\label{lem1}
\end{lem}

{\bf D\'emonstration}

Appliquer le th\'eor\`eme de  K.Saito  (1.9) \cite{S3}.  Les
vecteurs $ \xi_0, \cdots, \xi_{\mu-1} $ forment un
syst\`eme     involutif      sur     ${\mathcal     O}_{{\bf
C}^{\mu}},$ i.e. il existe $c_{i,j}^k(s) \in {\mathcal
O}_{{\bf C}^{\mu}}$ tel que
$$[\xi_i,\xi_j]=     \sum_{k=0}^{\mu-1} c_{i,j}^k(s) \xi_k.$$
Cette propri\'et\'e peut \^etre d\'eduit de la formule (2.5).
En plus, on a $D  =  \{\Delta_{\mu}(s)  =  \det   \;S(s)   =0\}.$
Donc  ces  vecteurs  satisfont  \`a la  condition   du   th\'eor\`eme
mentionn\'e ci-dessus. C.Q.F.D.

En suite nous consid\'erons l'int\'egrale

$$K^{\lambda}_{k,x^0}(s) = \int_{ Reg(\gamma_s)} {\bar x}^k (F({\bar x}+x^0,
s') + s_0)^\lambda dx, \leqno(2.2)'
$$
o\`u ${\bar x} = x-x^0.$ 

Afin d'\'etablir un syst\`eme analogue \`a $(2.3)_i$ etc. pour
$K_{k,x^0}^{\lambda}(s),$
nous reproduisons le raisonnement depuis la  Proposition ~\ref{prop1}.
Tout d'abord, nous avons les relations suivantes au lieu de $ (2.3)_i,$
$$
f_{0}(x^0,  s',s_0)\frac  {\partial}{\partial  s_0}K^{\lambda}_{k
+i,x^0}(s)+ $$
$$+ \sum_{\ell =1}^{\mu }
f_{\ell}(x^0, s') \frac{\partial}{\partial s_0}
K^{\lambda}_{k+ \ell +i, x^0}(s) +
\frac{\partial}{\partial s_0} K^{\lambda}_{\mu  +  k  +  i  +  1,
x^0}(s)
= {\lambda}K^{\lambda}_{k +i ,x^0}(s)  ,
\;\;\; 0 \leq i \leq \mu,\leqno(2.3)'_i$$
avec les polyn\^omes $f_\ell(x^0,s'): =\frac{1}{\ell !}
(\frac{\partial}{\partial x})^\ell F(x,s') \mid_{x= x^0}, \ell \geq 1, $
et $ f_0(x^0,s',s_0):= F(x,s')+s_0.  $

Au lieu de $(2.4)_j$ on a, 
$$\sum_{\ell =1}^{\mu } {\ell}f_{\ell}(x^0,s') \frac{\partial}{\partial s_0}
K^{\lambda}_{\ell+ k + j, x^0}(s) +(\mu +1) \frac{\partial}{\partial s_0}
K^\lambda_{\mu  +  k+  1+j,x^0}(s)  =  -(k  +  j+1)K^{\lambda}_{k
+j,x^0}(s)  ,
\;\;\; -1\leq j \leq \mu -1. \leqno(2.4)'_j$$
Par la suite une notation simple $f_\ell$ remplacera
$f_{\ell}(x^0, s') .$

Avant de formuler l'\'enonc\'e sur l'op\'erateur 
qui annule l'int\'egrale $K^{\lambda}_{k ,x^0} (s), $  on d\'efinit le
d\'eterminant d'une matrice

$${\bf   P}(s,\frac   {\partial}{\partial    s_0})
= \left[
\begin{array}{cccc}
p_{0,0}&p_{0,1} & \cdots & p_{0,\mu -1} \\
p_{1,0}& \vdots & \cdots & p_{1, \mu -1}\\
\vdots& \vdots & \cdots & \vdots  \\
p_{\mu -1,0}& p_{ \mu -1 ,1} & \cdots & p_{\mu -1,\mu -1} \\
\end{array}
\right],
$$ avec    les
composantes non-commutatives comme suit,
$$ det {\bf P} (s,\frac {\partial}{\partial s_0}):=$$
$$  = \sum_{i_0, i_1, \cdots, i_{\mu -1}}
sign (i_0,i_1, \cdots, i_{\mu -1})p_{i_{\mu -1}, \mu -1}\cdots  p_{i_1,1}
p_{i_0,0}, \leqno(2.6)$$
pour $ p_{i,j} = p_{i,j}^1(s)\frac {\partial}{\partial s_0}+ p_{i,j}^0(s), $
avec $ p_{i,j}^1(s),  p_{i,j}^0(s)  \in {\bf R}[s].$
Dans (2.6), l'indice $ (i_0,i_1, \cdots, i_{\mu -1}) $ parcourt toutes les
permutations de $(0,1, \cdots, \mu -1).$
Par la suite on se servira de la notation $ \frac
{\partial}{\partial s} = $
$(\frac     {\partial}{\partial     s_0},      \cdots,      \frac
{\partial}{\partial s_{\mu-1}}).$

\begin {prop}

i)

L'op\'erateur diff\'erentiel ${\mathcal P} (s,\frac {\partial}{\partial s})$
d'ordre $\mu$ d\'efini ci-dessous annule l'int\'egrale
$K_0^{\lambda}(s),$
$$ {\mathcal P} (s,\frac {\partial}{\partial s})K_0^{\lambda}(s) =0,$$
o\`u
$$ {\mathcal P} (s,\frac {\partial}{\partial s}): =
det  (S  \frac   {\partial}{\partial   s_0}   -   L   -V(s'))  $$
les matrices $S,L ,  V(s')$ sont celles d\'efinies dans la Proposition
~\ref{prop2}.

ii)
L'op\'erateur ${\mathcal P}^{(k)}_{x^0}
(s, \frac {\partial}{\partial s}) $ d'ordre $\mu +1$ qui  annule l'int\'egrale
$K^{\lambda}_{k,x^0}(s)$ de $(2.2)$
admet la repr\'esentation d'une mani\`ere analogue \`a
${\mathcal P} (s,\frac {\partial}{\partial s})$,

$$    {\mathcal    P}^{(k)}_{x^0}    (s,\frac     {\partial}{\partial
s})K^{\lambda}_{k,x^0}
(s)= 0,$$
$$ {\mathcal P}^{(k)}_{x^0} (s,\frac {\partial}{\partial s}) =
det ({\tilde S }\frac {\partial}{\partial s_0} - {\tilde L}
-{\tilde V(s')}) +
\sum_{j=1}^{\mu} {\tilde T'}_j (s,\frac {\partial}{\partial s})
,$$
o\`u les matrices ${\tilde S }, {\tilde  V(s')}$  sont  d\'efinies  \`a
partir de $ S,  V(s'):$
$$ {\tilde S } =$$
\begin{center}
$\left[\begin{array}{ccccc}
s_{0} -\tilde {s_0} & f_1(x,s') & \cdots & f_{\mu-1}(x,s')& f_{\mu} (x,s') \\
0& & &  &\\
\vdots &  & S'(s) & &  \\
0& &  & & \\
\end{array}
\right].$
\end{center}
$$ {\tilde L} = {\rm diag} (\lambda + \frac{k}{\mu +1},
\ldots, \lambda + \frac{k + \mu}{\mu +1}).$$
$$
{\tilde V(s') } =
\left[
\begin{array}{ccc}
0& \cdots &0 \\
 0 & \ddots    &\vdots \\
 V'(s')& 0 & 0 \\
\end{array}
\right],
 $$
$$ S'(s) = B_f \cdot S(s) \cdot C_f$$
avec $ B_f, C_f \in GL(\mu, {\bf R}),$ et $V'(s') \in  End  ({\bf
R}^{\mu -2}) \otimes {\bf R}[s'].$  Ici, on note ${\tilde s}_0 = -F(x^0,s').$
Les op\'erateurs ${\tilde T'}_j (s,\frac {\partial}{\partial s})$
d'ordre $\mu$ sont d\'efinis  dans $(2.9)'$ ci-dessous.
\label{prop3}
\end{prop}

\begin{remark}
{\em L'ordre de l'op\'erateur
${\mathcal    P}^{(k)}_{x^0}    (s,\frac     {\partial}{\partial
s}),$  \'egal \`a $\mu +1 ( >\mu),$ s'explique par l'existence d'un
cycle de plus $\gamma^0$ qui d\'efinit
$K^{\lambda}_{k, x^0}(s).$  Notons  $  \{  x^{(j)}\}  _{j  =  1,
\cdots, \mu +1}= \{ x \in {\bf C}: F(x,s') +s_0 =0 \}. $
Alors      $K_{0}^{\lambda}(s)$        a        ses        cycles
\'evanescents
$ \gamma^j = \{x^{(j)}  - x^{(j-1)} \} , j =2, \cdots, \mu +1.$ 
Par contre $K^{\lambda}_{k, x^0}(s)$
 peut \^etre d\'efini le long d'un cycle $x^0 - x^{(1)}$
\`a part des cycles \'evanescents associ\'es \`a  $K_{0}^{\lambda}(s).$
}
\label{remark2}
\end{remark}

{\bf D\'emonstration de la Proposition~\ref{prop3} }

{\bf Preuve de i)}
Nous notons la relation $(2.5)$ par
$${\bf     P}(s,\frac      {\partial}{\partial       s_0})    {\bf
K}^{\lambda}  (s)=  (S  \frac    {\partial}{\partial  s_0}  -L  -
V(s') ) {\bf K}^{\lambda}(s)=0. $$
Notons
$$ {\bf P} (s,\frac {\partial}{\partial s_0})= (p_{i,j}) _{0 \leq
i,j \leq \mu -1 } $$  avec les composantes $ p_{i,j} =
\sigma_{i,j}(s)\frac {\partial}{\partial s_0}+ p_{i,j}^0(s').$
Ici  les $ \sigma_{i,j}(s)$ sont introduits dans le lemme~\ref{lem1}.
Selon cette notation, (2.5) s'\'ecrit
$$ p_{j,0}K^{\lambda}_0(s) +  p_{j,1}K^{\lambda}_1 (s) +\cdots
+ p_{j,\mu -1} K^{\lambda}_{\mu -1} (s) =0, 0 \leq  j  \leq  \mu
-1. \leqno(2.7)_j $$
On prend une combinaison des expressions $(2.7)_j$,
$$ \sum_{j = 0}^{\mu -1 } sign (j, i_1, i_2,  \cdots,
i_{\mu-1})p_{i_{\mu-1},  \mu
-1} \cdots p_{i_1, 1}(2.7)_j. \leqno(2.8)$$
Alors il est possible de voir que l'op\'erateur devant $K^{\lambda}_0(s)$
doit \^etre \'egal \`a $det {\bf P}(s).$

D'autre part, l'op\'erateur devant $K^{\lambda}_k(s) $ est
$$ \sum sign (i_0, i_1, i_2,  \cdots,  i_{\mu-1})p_{i_{\mu-1},  \mu
-1} \cdots p_{i_k k} \cdots p_{i_1, 1} p_{i_0 k}=0.  \leqno(2.9)$$
Autrement dit, $(2.8)$ co\"incide avec
$$ det (S  \frac    {\partial}{\partial  s_0}  -L  -
V(s') ) K_0^{\lambda}(s).$$
Ainsi on a d\'emontr\'e i).

{\bf Preuve de ii) }

Afin de trouver un op\'erateur ${\mathcal P}^{(k)}_{x^0}(s,\frac
{\partial}{\partial s}) $
annulant $K^{\lambda}_{k, x^0}(s), $ on part de la relation
$\tilde {\bf \Sigma} \cdot  \tilde{b} = \tilde{a}$ avec
\begin{center}
$\tilde {\bf \Sigma}=
\tiny{ \left [
\begin{array}{cccccccccccc}
f_0&f_1&\cdots& f_{\mu -2} & f_{\mu-1}& f_{\mu}&1&0&\cdots&0&0&0\\
0&f_0&\cdots &f_{\mu-3}&f_{\mu -2} & f_{\mu-1}& f_\mu &1&\cdots&0&0&0\\
0&0&\cdots&f_{\mu -4} &f_{\mu-3}&f_{\mu -2} & f_{\mu-1}& f_\mu&\cdots&0&0&0\\
0&0 &\cdots&\ddots&\vdots &\vdots&\ddots &\ddots&\cdots& 1 &0 & 0 \\
\vdots&\vdots&\cdots&\ddots&\vdots&\vdots&\ddots  &\ddots&\cdots&
f_{\mu}&1 & 0 \\
0&0&\cdots&0& f_0& f_1&f_2&f_3& \cdots &f_{\mu-1}& f_\mu&1 \\
0&f_1&2f_2&\cdots&(\mu -1)f_{\mu-1} & \mu f_\mu &\mu +1&0&\cdots&0&0\\
0&0&f_1&\cdots&(\mu-2)f_{\mu-2}
&({\mu} -1)f_{\mu-1} & \mu f_\mu &{\mu +1}&0&\cdots&0&0\\
0&0 &&&\vdots &\vdots&\ddots &\ddots&\ddots& \ddots &0 & 0 \\
0&0 &\cdots
&0&f_1&2f_2&\cdots&(\mu-2)f_{\mu-2}&({\mu}-1)f_{\mu-1}& \mu
f_\mu&{\mu +1}&0\\
0&0 &\cdots&0
&0&f_1&\cdots&(\mu -3)f_{\mu -3}
&(\mu-2)f_{\mu-2}&({\mu}-1)f_{\mu-1}&\mu f_\mu&{\mu +1}\\
\end{array}
\right] } ,$
\end{center}
$$\tilde{a}= ^t(
{\lambda}K^{\lambda}
_{k, x^0},{\lambda} K^{\lambda}_{k+1, x^0 },\cdots,{\lambda}K^{\lambda}_{k+
\mu, x^0 },
-(k+1) K^{\lambda}_{k , x^0}, -(k+2)K^{\lambda}_{k+1 , x^0},\cdots, -(k+ \mu
+1)K^{\lambda}_{k + \mu, x^0 } )  $$
$$ \tilde{b}= ^t (\frac{\partial}{\partial s_0}K^{\lambda} _{k, x^0},
\frac{\partial}{\partial s_0}K^{\lambda}_{k+1, x^0},\cdots,
\frac{\partial}{\partial s_0}K^{\lambda}_{k+ \mu , x^0}).$$
Apr\`es des calculs analogues au cas
${\bf K}^{\lambda}$, pour  $\mu +1$ int\'egrales
${\bf K}^{\lambda}
_{k , x^0}:=$ $ ^t(K^{\lambda}
_{k, x^0},$ $ K^{\lambda}_{k+1, x^0},$ $\cdots,$ $
K^{\lambda}_{k+\mu , x^0}),$ on voit
$$({\tilde S} \frac {\partial}{\partial s_0}- {\tilde L'} -{\tilde
V}(s')) {\bf K}^{\lambda}_{k, x^0} =0 \leqno(2. 10)$$
pour
$$ {\tilde S } =$$
\begin{center}
$\left[\begin{array}{ccccc}
s_{0} -\tilde {s_0} & f_1(x,s') & \cdots & f_{\mu-1}(x,s')& f_{\mu} (x,s') \\
0& & &  &\\
\vdots &  & S'(s) & &  \\
0& &  & & \\
\end{array}
\right],$
\end{center}
avec $ S'(s) = B_f \cdot S(s)  \cdot  C_f$  o\`u  $B_f,  C_f  \in
GL(\mu, {\bf R})$ dont les composantes sont d\'etermin\'ees par les
coefficients de $f_k(x^0,s'), k = 1,\cdots, \mu.$
Les autres matrices sont donn\'ees comme suit:
$$ {\tilde L} = {\rm diag} (\lambda + \frac{k}{\mu +1},
\ldots, \lambda + \frac{k + \mu}{\mu +1}), $$
et
$$
{\tilde V(s') } =
\left[
\begin{array}{ccc}
0& \cdots &0 \\
 0 & \ddots    &\vdots \\
 V'(s')& 0 & 0 \\
\end{array}
\right],
 $$
avec $V'(s') \in  End  ({\bf
R}^{\mu -2}) \otimes {\bf R}[s'].$
De (2.10), on d\'eduit
$$ det ({\tilde S}  \frac    {\partial}{\partial  s_0}   -{\tilde
L}  - {\tilde V}(s') ) K_{k,x^0}^{\lambda}(s) +  \sum_{j =1}^{\mu}
{\tilde T}_j (s,\frac{\partial}{\partial
s_0})K^{\lambda}_{k +j, x^0 }  (s) =0 .$$

Ici, si on note $({\tilde S}  \frac    {\partial}{\partial  s_0}
-{\tilde L}  - {\tilde V}(s') ) = ( {\tilde p}_{i ,j  })_{0  \leq
i, j \leq \mu}, $ l'op\'erateur devant $K^{\lambda}_{k +j,  x^0  }
(s)$ est d\'efini comme suit:
$${\tilde T}_j (s,\frac {\partial}{\partial  s_0})
:= \sum sign (i_0, i_1, i_2,  \cdots,  i_{\mu}){\tilde p}_{i_{\mu},
\mu } \cdots {\tilde  p}_{i_1,  1}  {\tilde
p}_{i_0, j}. \leqno(2.9)' $$
L'existence     de l'op\'erateur     ${\tilde     T}'_j     (s,\frac
{\partial}{\partial  s})$ t.q.
$${\tilde T}'_j (s,\frac
{\partial}{\partial  s})K^{\lambda}_{k, x^0 }  (s)
=
{\tilde T}_j (s,\frac
{\partial}{\partial  s_0})K^{\lambda}_{k +j, x^0 }  (s),$$
d\'ecoule des relations
$$ (\frac {\partial}{\partial  s_0})K^{\lambda}_{k +j, x^0 }
(s)  = (\sum _{i=0}^j B_{ji}^{(k)}
\frac {\partial}{\partial  s_i}) K^{\lambda}_{k, x^0 } (s),
\;\; 1 \leq j \leq \mu -1 ,$$
o\`u les $B_{ji}^{(k)} $ sont d\'etermin\'es par
$$ \sum_{i=0} ^{j} B_{ji}^{(k)} ({\bar x} +x^0)^i {\bar x}^k
= {\bar x}^{k+j}.$$
Pour $ K^{\lambda}_{k + \mu, x^0 }, $
 $$ (\frac {\partial}{\partial  s_0})K^{\lambda}_{k + \mu, x^0 }
(s)  = [ -\sum _{j =0 }^{\mu -1} (B_{j}(x^0) + s_j)
\frac {\partial}{\partial  s_j} + \lambda ] K^{\lambda}_{k, x^0 }
(s),$$
o\`u les $B_{j}(x^0) $ sont t.q.
$$ ({\bar x} +x^0)^\mu -  {\bar x}^\mu
= \sum _{j =0 }^{\mu -1} B_{j}(x^0) ( {\bar x} + x^0)^j  .$$
C.Q.F.D.
 \vspace{2pc}

{\center{\section{\bf
Une version forte du th\'eor\`eme d'isomonodromie
}}}

Dans cette section, nous \'etablirons une version forte
du th\'eor\`eme d'isomonodromie. Il s'agit de l'invariance du comportement
asymptotique de  l'int\'egrale   $K^{\lambda}_{k}(s_0,s')$  (i.e.
celle  des expopsants caract\'eristiques) pr\`es
du point singulier (r\'egulier) $s_0= s_0(s')  $ lors de transition
le long d'une composante  stratifi\'ee de l'ensemble critique 
$D =\{ s \in {\bf C}^{\mu};\Delta_{\mu}(s) =0 \}.$ D'apr\`es le th\'eor\`eme
de Picard-Lefschetz et celui de Brieskorn \cite{Br}, la
monodromie du syst\`eme de Gauss-Manin (2.5) est identifi\'ee \`a
la monodromie locale
de Picard-Lefschetz:
$$ \Lambda_{\ast, s^\bullet}: H_1(X^{(\mu,\nu)}_s, {\bf C}) \rightarrow
 H_1(X^{(\mu,\nu)}_s, {\bf C}), \;\;\; s \not \in D, s^\bullet \in D,$$
qui se r\'ealise dans $GL(\mu,{\bf C}).$ Ici  on d\'enote par
$\Lambda_{\ast, s^\bullet}$ l'action du lacet $\Lambda \in
\pi_1({\bf C}^{\mu}\setminus D  )$ autour d'un point
$s^\bullet \in D$ qui agit sur le groupe d'homologie.

Avant de formuler le th\'eor\`eme, nous pr\'ecisons la notion des
exposants
caract\'eristiques.

\begin{dfn}
Soit $ (s',t(s'))  \in  D = \{  (s',t) \in {\bf C}^{\mu};  P_m(t,s')  =0
\},$  un point du lieu singulier $D$
d'une \'equation diff\'erentielle du type
de Fuchs (cf.\cite{AK} ):
$$P(s',t, \frac {\partial}{\partial  s'},
\frac {\partial}{\partial  t})=  [P_m(s', t)
(\frac {d}{d t})  + P_{m-1}(
s',t, \frac {\partial}{\partial  s'},
\frac {\partial}{\partial  t}) + \cdots + P_0(s', t)]u(s', t) =0,
\leqno(3.1)$$
o\`u
$$
P_{m-j}(s',t, \frac {\partial}{\partial  s'},
\frac  {\partial}{\partial   t})=    \sum   _{\alpha   \in   {\bf
N}^{\mu-1}, 0 \leq \mid \alpha \mid \leq m-j } P_{m-j, \alpha}
(s',t) ( \frac {\partial}{\partial  s'})^\alpha
(\frac {\partial}{\partial  t})^{m-j - \mid \alpha \mid}  $$
avec $\alpha = (\alpha_1, \cdots, \alpha _{\mu -1}),$ $ \mid  \alpha
\mid = \alpha _1 + \cdots + \alpha _{\mu -1}, $
$ (\frac {\partial}{\partial  s})^\alpha =
(\frac {\partial}{\partial  s_1})^{\alpha_1}  \cdots
(\frac {\partial}{\partial  s_{\mu -1}})^{\alpha_{\mu -1}}.$
L'op\'erateur du type de Fuchs
$P(s',t, \frac {\partial}{\partial  s'},
\frac  {\partial}{\partial   t})$ est  dit   \`a   multiplicit\'e
$\kappa$ le long de $D$ pr\`es  du point  $  (\dot{s}',t(\dot{s}'))
\in D,$ si  les  conditions   suivantes   sont   satisfaites
pour   ses
coefficients $P_{m-j, \alpha}(s',t) .$ On comprend par $t = t(s')$
l'\'equation locale de $D.$
D'abord on demande que la d\'ecomposition 
$$ P_{m}(s',t) = (t - t(s'))^{\kappa} Q_{ m }(s',t),$$
ait lieu pr\`es    de    $(\dot{s}',    t(\dot{s}'))$    avec    $Q_{    m
}(\dot{s}',t(\dot{s} ')) \not = 0.$
Deuxi\`emement, on suppose que $(t-t(s'))^{\kappa -j} \mid
P_{m-j, \alpha}(s', t)$ pour tout $0 \leq j \leq \kappa $ et
$\alpha \in {\bf N}^{\mu -1}.$
Alors il est possible de voir  que  l'expression  suivante  donne
naisssance \`a un polyn\^ome en $\rho,$
$$ \Pi_0(\rho, s')= (t- t(s'))^{-\rho}[ P(s',t, \frac
{\partial}{\partial  s'}, \frac {\partial}{\partial   t})
] (t- t(s'))^{\rho}\mid_{t = t(s')} $$
$$ = Q_{m}(s' ) \rho (\rho -1) \cdots (\rho -m+1) + Q_{m-1} (s') \rho (\rho
-1)\cdots(\rho -m+2) + \cdots+$$
$$+  Q_{m-  \kappa  +  1}   (s')  \rho  (\rho  -1)  \cdots  (\rho  -m+
\kappa) +  Q_{m-  \kappa }   (s')  \rho  (\rho  -1)  \cdots  (\rho  -m+
\kappa +1), \leqno(3.2)$$
pour une collection  de  fonctions  semi-alg\'ebriques  $Q_m(s'),
\cdots, Q_{m- \kappa}(s'). $

L'\'equation alg\'ebrique   en
$\rho, $ $ \Pi_0(\rho, s') =0$ se nomme
l'\'equation d\'eterminante relative  de l'op\'erateur
$P (s',$ $ t ,$ $\frac
{\partial}{\partial  s'},$ $\frac {\partial}{\partial   t})$
au point $(s', t(s')).$
\label{dfn31}
\end{dfn}

\begin{lem}
Au voisinage de son point singulier
r\'egulier $t = t(s') ,$ le comportement asymptotique des
m-solutions $u_1(s',t),
\cdots, u_m(s', t)$ de l'\'equation
du type de Fuchs (3.1) est d\'et\'ermin\'e
exactement (pas modulo ${\bf Z}$) par les exposants caract\'eristiques de
$P(s',t, \frac
{\partial}{\partial  s'}, \frac {\partial}{\partial   t})$
\`a $t =t(s').$
Cela veut dire que comme $m$ solutions on peut prendre
$$ u_{j,\ell}(s',t) \sim (t- t(s'))^{\rho_j}(ln (t- t(s')))^{\ell}
\sum_{k \geq 0}a_k^{(j,\ell)}(s') (t- t(s') )^k, \;\;\; 0 \leq \ell
\leq L_j-1,$$
o\`u $L_j$ est la multiplicit\'e de
racine $\rho_j$ dans l'\'equation d\'eterminante $(3.2),$ $ \sum_j L_j = m.$
En plus $a_0^{(j,\ell)}\not = 0.$
\label{lem31}
\end{lem}

Quant \`a d\'emonstration de cet \'enonc\'e, ce n'est qu'une
modification \`a la version en plusieurs variables (une variable
plus  plusieurs  param\`etres)  des  r\'esultats  classiques   sur
l'\'equation du type de  Fuchs.  A  ce  propos,  je  renvoie  les
lecteurs au livre classique de Coddington-Levinson, ou d'Ince
qui explique la m\'ethode de Frobenius.

Avant de formuler le  th\'eor\`eme d'isomonodromie, nous
\'etablirons  ici
un  lemme  et  y  introduisons  la  notion  de   stratification
logarithmique.
 Notons par $Der_ { {\bf C}^{\mu}}(\log D)(p), p\in D $ le sous-espace
vectoriel de  l'espace  tangent  $T_p   {\bf  C}^{\mu}   $   qui
consiste des
valeurs $\delta(p)$ de vecteurs logarithmiques $\delta \in Der_ { {\bf
C}^{\mu}}(\log D),$ introduits par K.Saito \cite{S3}. Ici $\delta$ est
un  vecteur
\`a coefficients polynomiaux tangent \`a $D.$

D'abord  on d\'efinit l'ensemble $D ( i_1 ,$ $\cdots, $ $ i_{k+1})$ comme
un sous-ensemble de $D$ sur lequel les $(k+1) -$ vecteurs
logarithmiques $ \xi_{i_1},$ $\cdots,$ $\xi_{ i_{k+1}} $ $\in$ $ Der_ { 
{\bf C}^{\mu}}(\log D) $ sont exprim\'es par
$ \xi_{j_1},$ $\cdots,$ $\xi_{ j_{\mu -k -1}}$ avec
$ \{ j_1,$ $ \cdots,$ $ j_{\mu - k-1} \}$
$ = \{0, $ $\cdots, $ $\mu-1 \}$  $ \setminus \{ i_1, $ $\cdots, $ $i_{k + 1}
\}.  $
C'est  \`a  dire,  si  $s  \in  D(i_1,\cdots,   i_{k+1}),$  alors  il
existe    des    fonctions    semi-alg\'ebriques non nulles
$A^{p}_{j_{\ell}}(s'), $$A^{p}_{i_p}(s') $
$\in {\bf C} [s_0(s'),s'], $ $ 1 \leq $ $p \leq $ $k +1,$ $ 1 \leq  $  $\ell  \leq
$ $\mu -k -1$ t.q.
$$ A^{p}_{i_p}(s')\xi_{i_p}(s) = \sum_{\ell  =1}^{\mu  -k  -1}  A^{p}_{j_{\ell}}(s')
\xi_{j_{\ell}}(s), 1 \leq p \leq k+1. $$
sur $\pi (D(i_1,\cdots,   i_{k+1})).$ Ici $\pi:{\bf C}^{\mu}_s \rightarrow
{\bf C}^{\mu-1}_{s'} $ d\'enote la projection sur ${\bf C}^{\mu-1}_{s'}$.
En fait $A^{p}_{\bullet}(s')$ sont donn\'ees par des mineurs $ (\mu  -k  -1)
\times (\mu  -k  -1)$ de $S(s_0(s'),s').$

 On se place dor\'enavant dans le  compl\'ement de 
l'ensemble de Maxwell du polyn\^ome
(2.1). L'ensemble de Maxwell $M \subset D$ est l'ensemble des  $s$
pour lesquels
l'\'equation $F(z,s') + s_0 =0$ poss\`ede plusieurs points critiques
ind\'ependants qui donnent
la m\^eme valeur critique (ex.$ \mu =3, s_0= \frac{s_2^2}{4}, s_1=0$  ).
Notons
l'ensemble de Maxwell par $M$ et $D_M: = D \setminus M.$ L'ensemble de Maxwell
lui m\^eme est stratifi\'e.

\begin{lem}

i)   Le point $p$ appartient \`a l'ensemble  $D_{M}(i_1,$  $\cdots,$  $
i_{k+1})$
si et seulement s'il existe $\mu -k -1$  vecteurs
lin\'eairement ind\'ependants $\delta_{j_1}(p), \cdots,
\delta_{j_{\mu -k -1}}(p) \in  Der_ { {\bf
C}^{\mu}}(\log D)(p). $
   Cela veut dire que l'ensemble $D_{M}(i_1,\cdots,  i_{k+1})$ est
d\'etermin\'e uniquement par le nombre $\mu - k -1$ des vecteurs
logarithmiques ind\'ependants sur
lui. On note tel ensemble  par ${D_M}^{(k)}.$

ii)  On a la structure de stratification  comme suit:
$$  D_M= \bar{ D_M}^{(0)}  \supset \bar  {D_M}^{(1)}  \supset
\bar{ D_M}^{(2)}  \supset  \cdots \supset D^{(\mu -1)} =\{ 0\}.$$
L'ensemble $D_M^{(k)}$ est  une strate    de  codimension  un  dans
$\bar{ D_M}^{(k-1)}.$  A chaque point $p$ de la  strate  $D_M^{(k)},$  on  a
$rang
Der(  log   D)(p)  =  \mu  -  k   -1.$   ({\it   La   stratification
logarithmique} de K.Saito, \cite{S3}.)

iii) La strate $D_M^{(k)}$ peut \^etre d\'efini par les 
mineurs $(\mu
-k) \times (\mu - k)$ de la matrice $S(s).$

iv) Si $s \in D_M^{(k)},$ alors $F(z,s') + s_0 =0$  poss\`ede  une
racine d'ordre
$k+2$.
\label{lem2}
\end{lem}

\begin{remark}
Si on prend une strate $M^{\alpha}$ de l'ensemble de Maxwell et un point
sur lui $p \in M^{\alpha},$ alors  $rang
Der(  log   D)(p)   > \mu -1 - codim_D M^{\alpha}.$
\end{remark}

{\bf D\'emonstration}

   Pour d\'emontrer l'\'enonc\'e i), on proc\`ede par 
induction.  Tout d'abord, par d\'efinition de
$D (i_1,$  $\cdots,$  $i_{k+1})$ il est \'evident que
$D_{M}(\mu -1)$ $=$ $D_{M}(\ell),$ $\ell = 0, \cdots, \mu -2.$
Ensuite supposons
que  $s  \in
D_M(\mu -2, \mu -1 ) \subset D_M$ et montrons que
$s \in D_{M}(i_0,$  $i_{1}),$ pour n'importe quelle paire
$ \{ i_0, i_{1} \}$ $ \subset \{0, \cdots, \mu-1 \}.$
L'hypoth\`ese $s \in D_M(\mu -2, \mu -1 )$ entra\^{i}ne qu'il existe
deux
collections de fonctions semi-alg\'ebriques non nulles $(c_{0}(s'),   \cdots, c_{\mu -1}(s'))$
$\in {\bf C} [s_0(s'),s']$
 et $(c_{0}' (s') ,$ $ \cdots,$ 
$ c_{\mu -2}' (s'),0 )$$\in {\bf C} [s_0(s'),s']$
telles que
$$  (c_{0}(s'),    \cdots,   c_{\mu   -1}(s'))   \cdot   S(s_0(s'),s')   =0
\leqno (3.3)$$
$$ (c_{0}' (s'),    \cdots,  c_{\mu  -2}'(s'),0)  \cdot  S(s_0(s'),s')  =0,
\leqno (3.3)'$$
pour $s' \in \pi (D_M(\mu -2, \mu -1 )  ). $
Et il est impossible de trouver de fonctions semi-alg\'ebriques non nulles
$ (c_{0}''(s'), \cdots , \check{0}^{ j_1}, \cdots,
\check{0}^{ j_2 } , \cdots, c_{\mu -1}''(s'))  $ telles que
$$(c_{0}''(s'), \cdots , \check{0}^{ j_1} , \cdots,
\check{0}^{ j_2 } , \cdots, c_{\mu -1}''(s')) \cdot  S(s_0(s'),s') =0, \leqno(3.4)$$
pour $s' \in \pi (D_M(\mu -2, \mu -1 )  ). $
Il  faut  remarquer  que   la   derni\`ere   condition   entra\^{i}ne
qu'aucun de  mineurs $(\mu -2) \times  (\mu -2) $
de $S(s)$ ne s'annule. Car si pour
les indices $ \{ k_0, \cdots, k_{\mu  -3} \},$ $ 0  \leq  k_0
<  \cdots $ $<  k_{\mu -3} \leq \mu -1,$
les 
mineurs $(\mu -2) \times  (\mu -2) $ s'annulent i.e.
$$ S(^{k_0, \cdots,  k_{\mu -3}}_ { j_3, \cdots, j_{\mu}})=0,$$
avec
 $ \{ j_3, \cdots, j_{\mu} \}$
$ = \{1, \cdots, \mu \}$  $ \setminus \{ j_1, j_2 \},$ alors on peut
trouver une collection de fonctions comme $(3.4).$

   Revenons  \`a  la  condition  de  d\'eg\'ener\'escence   de  $\mu  -1$
vecteurs logarithmiques.  Pour  obtenir  une  relation  lin\'eaire
non-triviale entre $\xi_0, \cdots, \xi_{i_0 -1}, \xi_{i_0  +1},
\xi_{\mu -1}, $ il est suffisant de prendre la diff\'erence entre
$ c_{i_0}' (3.3)$ et $ c_{i_0} (3.3)'. $ Il est facile de  se
convaincre que  tous  les  coefficients
$c_{i_0}'  c_0  - c_0' c_{i_0},$ $ \cdots,$
$c_{i_0}'  c_{\mu -2}  - c_{\mu -2}'c_{i_0},$
$ c_{\mu -1} c'_{i_0} $ sont non nuls sauf celui correspondant \`a
$\xi_{i_0}.$ Car, sinon, il existe un 
mineur $(\mu -2) \times  (\mu -2) $
d\'eg\'en\'er\'e de $S(s)$ et on est alors dans 
la situation de $(3.4).$ 
On a ainsi obtenu l'expression suivante avec certaines fonctions
$A^{1}_{i_1}(s'), A^{1}_{j_{\ell}}(s') \not = 0,$
$$A^{1}_{i_1}(s') \xi_{i_1}(s_0(s'), s') = \sum_{j_\ell  \not = i_0, i_1} 
A^{1}_{j_{\ell}}(s')\xi_{j_{\ell}}(s_0(s'), s'), $$
pour $i_1 \not = i_0.$
Quant \`a  l'expression semblable pour $\xi_{i_0},$ on prend
la diff\'erence entre
$ c_{i_1}' (3.3)$ et $ c_{i_1} (3.3)'. $ Ainsi l'\'enonc\'e $i)$ est
d\'emontr\'e pour $k= 1.$

 Les d\'emarches d\'emand\'ees pour accomplir la  d\'emonstration
de $i)$ pour la strate de codimension $k >1$ 
sont analogues \`a l'argument ci-dessus.

Bien \'evidemment, cette d\'emarche ne s'applique pas aux strates de
Maxwell. Par
example la strate g\'en\'erique $M^{(0)}$ de  $M$  est de  codimension
un dans $D,$
mais il existe $\mu -1$ vecteurs logarithmiques
$\xi_1, \cdots, \xi_{\mu -1}, $ qui sont lin\'eairement ind\'ependants
sur   $M^{(0)}.$

L'\'enonc\'e $ii)$ est un corollaire de $iii).$

La d\'emonstration de l'\'enonc\'e $iii)$  s'appuie  sur  l'argument
suivant.  Si  un  des   mineurs $k  \times  k  $  de  la  matrice
$S(s)$ s'annulent 
au point $s \in D  $  alors  tous  les  autres 
mineurs
s'annulent  aussi \`a ce point l\`a. Par exemple, m\^eme si
la strate $D (i_0, i_{1}) $ est d\'efinie par $  _{\mu  -1}  C_
{\mu -2} $ mineurs dans $D,$ la codimension de celui-ci dans  $D$
est \'egale \`a 1. Cela appara\^{i}t dans la d\'emonstration de $i).$
Ainsi la strate $D (i_0, \cdots, i_{k}, i_{k+1}) $
est un ensemble de  la codimension un
 de la  strate  $D (i_0, \cdots, i_{k}). $

Par cons\'equence de $ iii), $ l'\'equation $\Delta_\mu(s_0,s')=
det \;S(s)=0$ poss\`ede une racine d'ordre $(k+1)$ $s_0= s_0(s')$ sur
$D^{(k)}_M,$ ce qui donne 
$iv)$ par les r\'esultats classiques sur le  discriminant. C.Q.F.D.

Maintenant nous sommes susceptibles de  formuler  le  th\'eor\`eme
d'isomonodoromie renforc\'e pour les op\'erateurs introduits 
dans la Proposition
~\ref{prop3},
$$ {\mathcal  P}(s, \frac {\partial}{\partial s})=
\Delta_{\mu}(s) (\frac{\partial}{\partial s_0} )^{\mu}+
P_{\mu -1}(s, \frac{\partial}{\partial s} ) + \cdots
+ P_{1}(s, \frac{\partial}{\partial s} ) +  P_0(s)
$$
$$ =
\Delta_{\mu}(s) (\frac{\partial}{\partial s_0} )^{\mu}
+ \sum_{1 \leq j \leq \mu }
P_{\mu -j}(s) 
(\frac {\partial}{\partial  s_0})^{\mu -j},  $$
et
$$ {\mathcal  P}^{(\ell) }_{x^0}(s, \frac {\partial}{\partial s})
= (s_0 -{\tilde s}_0) \Delta_{\mu}
(s) (\frac{\partial}{\partial s_0} )^{\mu+1} +
P_{\mu, x^0}^{(\ell )}(s, \frac{\partial}{\partial s} ) + \cdots
+ P_{1, x^0}^{(\ell)} (s, \frac{\partial}{\partial s} ) +
P_{0,x^0}^{(\ell)}(s), $$
$$ =
(s_0  -{\tilde  s}_0)  \Delta_{\mu}(s)  (\frac{\partial}{\partial
s_0} )^{\mu +1}
+ \sum_{1 \leq j \leq \mu +1  }\sum   _{\alpha   \in   {\bf
N}^{\mu-1},} P_{\mu -j +1,\alpha , x^0}^{(\ell )}
(s) ( \frac {\partial}{\partial  s'})^\alpha
(\frac {\partial}{\partial  s_0})^{\mu +1  -j - \mid \alpha \mid},  $$
pour $\ell \geq 1.$ 

Introduisons une structure stratifi\'ee du  lieu de ramification
associ\'e
 \`a l'op\'erateur
${\mathcal  P}^{(\ell) }_{x^0}(s, \frac {\partial}{\partial s}).$ Notamment
nous d\'efinissons les strates,
$$ {\tilde D_M}^{(k)}:=
D_M^{(k)} \cap \{(s_0,s');s_0 = {\tilde s}_0 \}, \; k=-1,0,1
,...$$
$$ D_{M,-}^{(k)}:= D_M^{(k)}\setminus \{s_0 = {\tilde s}_0 \}.$$
 Eventuellement,
$$ {\tilde D_M}^{(-1)}:= \{s_0 = {\tilde s}_0 \} \setminus D. $$
Ici et dans la suite on note $\tilde s_0 = -F(x^0,s').$
\begin{thm} (Th\'eor\`eme d'isomonodromie renforc\'e)

i) L'\'equation d\'eterminante de  l'op\'erateur  ${\mathcal  P}(s,\frac
{\partial}{\partial s})$ est de  multiplicit\'e  constante  le
long de la  strate $D_M^{(k)}.$
C'est \`a dire au voisinage de chaque point singulier
$ (s_0(s'),s') \in D_M^{(k)}$
$$ P_{\mu}(s) = \Delta_{\mu}(s) = (s_0 - s_0(s'))^{k+1} Q_{\mu ,k}(s_0,s')$$
$$  P_{\mu  -j}(s)=   (s_0   -   s_0(s'))^{k   -j   +1}
Q_{\mu   -j,k}(s_0,s'), \;  1  \leq  j  \leq  k+1, $$
pour  $s' \in \pi ( D_M^{(k)} ).$
Ici les polyn\^omes $ Q_{\mu  -j, k} ( s_0 , s')$ sont des 
polyn\^omes  de
degr\'e $\mu -k -1 $ en variable $s_0$ t.q. $Q_{\mu,
k} ( s_0 , s')$ $ \not =$ $ 0 $ sur la strate $D_M^{(k)}.$

La m\^eme d\'ecomposition a lieu  au point
$  (\tau^{\mu+1}s_0,\tau^{\mu}s_1,
\cdots, \tau^2 s_{\mu -1}), \tau > {\bf R}_+.  $

ii) Pour l'op\'erateur ${\mathcal  P}^{(\ell) }_{x^0}
(s, \frac {\partial}{\partial s}),$ sur la strate ${\tilde D_M}^{(k)},$ on a
la factorisation  suivante de ses coefficients:
$$ P_{\mu +1 ,x^0 }^{(\ell)} (s)=
(s_0 -{\tilde s}_0)\Delta_{\mu}(s)= (s_0 -{\tilde s}_0)^{k+2}
Q_{\mu +1 ,k, x^0}^{\ell} (s_0,s')$$
$$  P_{\mu  -j  +1,\alpha, x^0  }^{(\ell) }(s)=  (s_0  -  {\tilde
s}_0)^{k   -j   +2}
Q_{\mu +1 -j, k, \alpha, x^0}^{\ell}(s_0,s'), \; 1  \leq  j  \leq
k+2,  \mid \alpha \mid = 0,1, $$
avec $Q_{\mu +1 ,k, x^0}^{\ell} (s_0,s') \not =0. $

Sur la strate $D_{M,-}^{(k)},$ on a la factorisation suivante pour les
coefficients:
$$  P_{\mu  -j  +1,\alpha, x^0  }^{(\ell) }(s)=
(s_0  -  s_0(s'))^{k   -j+1}
Q_{\mu +1 -j, k,\alpha, x^0}^{\ell}(s_0,s'), \; 1 \leq j \leq k+1,
\mid \alpha \mid = 0,1, $$
avec $s_0(s') \not = {\tilde s}_0$ et
$Q_{\mu +1 ,k, x^0}^{\ell}
(s_0,s') \not =0. $

iii) Les exposants caract\'eristiques de l'op\'erateur
${\mathcal  P}(s,\frac
{\partial}{\partial s})$
au point singulier
$ (s_0(s'),s')$ $\in$ ${ D_M}^{(k)}$  $\subset$ $ {D_M} $ ne
changent pas lors de la translation du point
$ (s_0(s'),s')$ vers un autre point
$ (s_0(t'),t')\in{ D_M}^{(k)}  \subset {D_M} $
 le long  de la  strate $  { D_M}^{(k)}  \subset {D_M}. $
D'une fa\c{c}on analogue,  les exposants caract\'eristiques de ${\mathcal
P}^{(\ell)}_{x^0}(s, \frac {\partial}{\partial s})$
ne changent pas lors de la translation le long des strates
${\tilde D_M}^{(k)},D_{M,-}^{(k)}   \subset {D_M} 
.$
\label{thm3}
\end{thm}

\begin{remark}

{\em C'est un \'enonc\'e  renforc\'e  d'une  assertion  du  type
suivant:}
``La repr\'esentation du groupe de monodromie autour du point $x$
dans $GL(\mu, {\bf C})$
ne  change  pas  lors  de  la translation  le  long  d'une  strate   de
l'ensemble critique
sur laquelle se trouve le point de d\'epart $x$.''

{\em
Cette assertion d\'ecoule   d'un  th\'eor\`eme
bien connu de E.Brieskorn (Satz 1, III , \cite{Br})}
"La  monodromie  de  l'op\'erateur   diff\'erentiel   singulier
$\nabla_{f,x}$ peut \^etre identifi\'e avec  la
monodromie
locale de Picard-Lefschetz de $f(z)=F(z,s')$ ( chez nous $F(z,s') =
z^{\mu+ 1} + s_{\mu-1} z^{\mu -1} + \ldots + s_1 z) $ \`a $x = -s_0.$ "
{\em La version  de   Picard-Lefschetz par F.Pham  pour
l'int\'egrand   ramifi\'e (\cite{Ph}, \cite{Vas}) permet
d'identifier   la   monodromie    de
Picard-Lefschetz gen\'eralis\'ee   avec    la    monodromie    de
l'op\'erateur diff\'erentiel qui annule les int\'egrales
correspondantes (i.e. chez nous $S \frac {\partial}{\partial
s_0}-L -V(s')$).}
{\em L'invariance des exposants de monodromie le long de la strate
$\mu = const.$ est connue depuis Varchenko \cite{Var2}. 
Pour autant qu'il s'agisse de 
l'int\'egrale hyperelliptique g\'en\'eralis\'ee, notre th\'eor\`eme est une
version forte de celui de Varchenko, car nous constatons l'invariance de 
$\mu$ (ou bien $\mu +1$ ) exposants pour chaque int\'egrale 
$  I_{\omega}   (t),$ pourtant Varchenko a d\'emontr\'e l'invariance
du minimum de ces $\mu$ exposants. 
}
\label{remark6}
\end{remark}

{\bf D\'emonstration du Th\'eor\`eme ~\ref{thm3}}

D'abord on montre l'\'enonc\'e i).

Nous commen\c{c}ons par l'\'enonc\'e sur la multiplicit\'e de
$P_{\mu  -j}(s).$
L'enonc\'e  sur  la  multitiplicit\'e  de   $P_{\mu}(s)$ est   une
cons\'equence   immediate   de la Proposition   ~\ref{prop3}   et
du Lemme~\ref{lem2}.
Pour d\'emontrer l'\'enonc\'e pour $P_{\mu  -j}(s), j >0,$
nous modifions la matrice
${\bf P} (s, \frac {\partial}{\partial s_0})
=  (S \frac {\partial}{\partial s_0} -L -V(s'))$
en une autre matrice ${\bf P^I}(s, \frac {\partial}{\partial
s_0}), {\bf I}= (i_1, \cdots,  i_j)
\subset \{ 0,1, \cdots, \mu -1 \}$ dont le $i_{\gamma}- $\`eme
rayon   est   remplac\'e   par   $   ( 0 ,   0 , \cdots ,    0 , \check{
1}^{i_{\gamma}} $ $, 0 , $ $\cdots , $ $0 ) , $ $\gamma  $  $=  $
$1, $ $\cdots, $ $j.$
D'apr\`es la Proposition ~\ref{prop3} i) on voit que
$P_{\mu  -j,  0}(s,  \frac  {\partial}{\partial  s_0})$  est  une
somme de
termes d'ordre $\mu -j$ des op\'erateurs $det
{\bf P^I}(s, \frac {\partial}{\partial s_0}),$ ${\bf I}$ $\in$ 
$\{$ $\mid {\bf I} \mid = j;
{\bf I} \subset \{ 0,1, \cdots, \mu -1 \} \}.$
Cet argument entra\^{i}ne que le coefficient de
$P_{\mu  -j}(s,  \frac  {\partial}{\partial  s_0})$
peut \^etre exprim\'e par une somme de  
mineurs $(\mu  -j  )\times   (\mu
-j)  $ de la matrice $S(s).$

Puisque $  S(s)  =  s_0  id  _{\mu}
+C(s'),$ et $rank S(s) = \mu  -  1-k$  pour  $s\in  D_M^{(k)},$  la
matrice $S(s)$ est diagonalisable pour  $s$  assez  proche  de  $
D_M^{(k)}.$
L'absence de blocs de Jordan non-triviaux est garanti par le Lemme~\ref{lem1}.
Autrement dit,  il existe une matrice  $B(s')  \in  GL(\mu,  {\bf
C}[s']) $ sur $ \pi(D_M^{(k)})$ telle que,
$$ B(s')^{-1} S(s) B(s') = diag (s_0 - s_{0,1}(s'), \cdots,  s_0 -
s_{0,1}(s'),
s_0 - s_{0,2}(s'), \cdots, s_0 - s_{0, \mu -k -1}(s')) $$
o\`u  $s_{0,i}(s') \neq s_{0,j}(s')$ si $i \neq j,$ le terme $s_0 -
s_{0,1}(s')$ appara\^{i}t $k +1$ fois.

Evidemment, les 
mineurs $(\mu -j )\times  (\mu  -j  )$ de la matrice
$  B(s')^{-1}  S(s)  B(s')$   contiennent   le   facteur   $(s_0   -
s_{0,1}(s'))^{k-j +1}, 0 \leq j \leq k+1.$ Ici on se souvient de
la formule de Binet-Cauchy qui exprime les  mineurs
$p \times  p$
d'une matrice $AB$ au  moyen  des  mineurs  $p  \times  p$ des
matrices $A$ et $B,$
$$ S(^{i_1, \cdots,  i_{p}}_ { j_1, \cdots, j_{p}})=
 \sum_{1 \leq k_1 < \cdots < k_{p}\leq \mu,
1 \leq \ell_1 < \cdots < \ell_{p}\leq \mu } B(^{i_1, \cdots,  i_{p}}_ {
     k_1, \cdots, k_{p}}) (B^{-1}S B ) (^{k_1, \cdots,  k_{p}}_ {
\ell_1, \cdots, \ell_{p}}) (B^{-1} ) (^{\ell_1, \cdots, \ell_{p}}
_{j_1, \cdots,  j_{p}}).$$
On en d\'eduit que les mineurs $(\mu -j) \times (\mu  -j)$ de la
matrice $  S(s) $ aussi  contiennent   le   facteur   $(s_0   -
s_{0,1}(s'))^{k-j +1}.$  Ainsi on a d\'emontr\'e
$$ P_{\mu -j}(s) = (s_0  -  s_{0,1}(s'))^{k-j  +1}(\mbox{polyn\^ome
de degr\'e } \mu -k -1\; \mbox {en }\; s_0), 0\leq j \leq k+1. $$

Ceci ach\`eve la d\'emonstration de $i)$ pour $\ell = 0.$

{\bf D\'emonstration de $iii)$ pour $\ell = 0.$}
Puisque la strate $D_M^{(k)}$ constitue en plusiers 
composantes disjointes 
la d\'emonstration sera divis\'ee en deux \'etapes.
Chaque   composante   est    r\'etractible    par    l'action    de
quasihomog\'enei\'et\'e,
$$ \tau : (s_0, \cdots, s_{\mu -1}) \rightarrow
( t^{\mu  +1}  s_0,  t^{\mu   }
s_1,\cdots,   t^2   s_{\mu -1}),  \;\; t \in {\bf C}^{\times}. \leqno(3.5) $$

Etape 1. Invariance sur une composante.

Supposons que les deux points $ (s_{0,1}(s'),s') ,  (s_{0,1}  (t'),
t') \in D_M^{(k)} $ \`a petite distance $\epsilon$ l'un de l'autre.

Gr\^ace   \`a   l'\'enonc\'e   $i),$   on   a   les   \'equations
d\'eterminantes   correspondant   \`a    ces    points    (voir
(3.2) ):
$$
\Pi_0(\rho,s') = \rho (\rho -1)  \cdots  (\rho  -  \mu  +  k  +2)
Q(\rho, s'),\leqno(3.6)$$
$$
\Pi_0(\rho,t') = \rho (\rho -1) \cdots (\rho - \mu + k +2)
Q(\rho,t'), \leqno(3.7)$$
o\`u $Q(\rho,s')$ un polyn\^ome de degr\'e $(k+1)$ en $\rho.$
Il est \'evident que les deux \'equations  $(3.6)$  et  $(3.7)$
poss\`edent les racines communes $\rho = \{ 0, 1, 2, \cdots,  \mu
-k -2\}.$ Par contre, chacune  a  des racines $\{ \rho_{\mu -k -1}(s')
\leq \cdots \leq \rho_{\mu -1}(s') \} $ et
$\{ \rho_{\mu -k -1}(t')\leq \cdots \leq \rho_{\mu -1}(t') \} $
qui peuvent \^etre diff\'erentes, mais pour certain petit
$\delta (\epsilon) >0 ,$
$$ \mid \rho_j (t')  -  \rho_j  (s')  \mid  <  \delta  (\epsilon).
\leqno(3.8) $$
D'apr\`es  le  th\'eor\`eme  de  E.Brieskorn - F.Pham  cit\'e   ci-dessus
(Remarque ~\ref{remark6}), on sait que
$$ e^{i\rho_{\mu -k -1}(s')}= e^{i\rho_{\mu -k -1}(t')},
\cdots, e^{i \rho_{\mu -1}(s')}= e^{i \rho_{\mu -1}(t')} .$$
\c{C}a veut dire que le racine $\rho_j(s')$ doit \^etre \`a
distance enti\`ere  
de $\rho_j(t')$.  La continuit\'e $(3.8)$ entra\^{i}ne
l'invariance des racines eux-m\^emes lors  de la transition
$$ (s_0(s'),s') \rightarrow (s_0 (t'), t') \in D_M^{(k)},$$ i.e.
$$ \rho_{\mu -k -1}(s')= {\rho_{\mu -k -1}(t')},
\cdots, {\rho_{\mu -1}(s')}= { \rho_{\mu -1}(t')} .$$

Etape 2. Aux points en sym\'etrie.

L'appartenance d'un point $(s',s_0(s'))$ \`a $D^{(k)}_M$ veut dire
qu'il existe un
polyn\^ome $q(z- \dot{z})$ de degr\'e $ \mu -k- 1,$ tel que $ q(0)\not = 0$
et $$ F(z,s') +s_0(s')=(z- \dot{z})^{k+2}q(z- \dot{z}),$$
pour certain $ \dot{z} \in {\bf C}. $
Gr\^ace \`a la r\'etraction $(3.5)$ et l'\'etape 1, il suffit de v\'erifier
la co\"incidence des exposants caract\'eristiques \`a des points avec
$s_0=1$ se trouvant sur les diverses composantes.

Notons $  \dot{z}_u = -\omega^{u},\;\;$
$ 0 \leq u
\leq k+1  $ pour $\omega = e^{i\frac{2 \pi i}{k+2}}.$
 On remarque qu'il existe deux \'ecritures  pour $F(z,s') + 1$ \`a ces points
i.e.
\begin{flushleft}
$F(z,s') + 1 = (z- \dot{z}_u)^{k+2}(q_{\mu -k-1,u}(z -\dot{z}_u)^{\mu -k-1}+
\cdots + q_{1,u}(z -\dot{z}_u) + q_{0,u})$
$=z^{\mu +1}+ s_{\mu -1,u}z^{\mu -1}+ \cdots + s_{1,u}z +1$
\end{flushleft}
Si on regarde les d\'eriv\'ees $(\frac{\partial}{\partial z})^i F(z,s'), \;0
\leq i \leq \mu+1  $ et si on compare les r\'esultats de deux expressions ci-dessus,
on obtient des \'equations lin\'eaires entre $(q_{0,u}, \cdots,q_{\mu
-k-1,u} )$
et $(s_{1,u}, \cdots, s_{\mu -1,u}).$ En r\'esolvant ces \'equations, on
obtient   $s_{j,u} = \omega^{u(\mu -j +1)}s_{j,0}$ et $q_{\mu -k -j,u}
= \omega^{u(j -1)}q_{\mu -k -j,0}.$
Cela veut dire que $D^{(k)}$ poss\`ede $k+2$ composantes
connexes qui sont repr\'esent\'ees par
$$ d_u = (1,s_{1,u}, \cdots, s_{\mu -1,u})= (1,\omega^{u \mu}s_{1,0},
\omega^{u(\mu -1)}s_{2,0},\cdots,  \omega^{2u}s_{\mu -1,0}), \;\;
0 \leq u \leq k+1.$$
Parmi les exposants en question, uniquement les $k+1$ exposants
non-entiers importent car les autres
sont forc\'ement
$\{0,1,\cdots, \mu -k -2\}  $ pour tous les $d_u.$ Donc, en vertu  de
la Proposition
~\ref{prop3}, afin de comparer les exposants aux points $d_u, 0 \leq u \leq
k+1,$ il suffit de comparer $ \rho^{(i)}_u$  des int\'egrales  ci-dessous
$$ I^{(i)}_u(t) = \int_{Reg(\delta_{i,u}(t))}(F(z,s') + 1-t)^{\lambda}dz \sim
v^{(i)}_{0,u} t^{ \rho^{(i)}_u  }( 1 + O( t ^{ \frac{1}{k+2}}) ).
$$
Ici $\delta_{i,u}(t)$ le cycle \'evanescent qui converge vers $\dot{z}_u $
lorsque $t \rightarrow 0, 1 \leq i \leq k+1 .$ On peut supposer que $
\delta_{i,u}(t)= \omega^u \cdot \delta_{i,0}(t). $
Pour analyser les exposants  $ \rho^{(i)}_u,$ faisons le changement de
variable analytique $z \rightarrow Z_u$ pr\`es de $z = \dot{z}_u$
$$ Z_u = w_{1,u}(z- \dot{z}_u) + \sum_{ i \geq 2} w_{i,u}(z- \dot{z}_u)^i$$
tel que $ F(z,s') + 1 = Z_u^{ k+2} . $ Ceci est possible vu que
$$   \frac{1}{(k+2)!}(\frac{\partial}{\partial z})^{k+2 -j} F(\dot{z}_u ,s')
= \delta_{0,j} q_{0,u}. $$ 
Alors
$$ I^{(i)}_u(t) = \int_{Reg(\tilde \delta_{i,u}(t))}(Z_u^{k+2}-t)^{\lambda}
(\frac{\partial Z_u}{\partial z} )^{-1}dZ_u \sim
 \sum_{j \geq 0}v^{(i)}_{j,u} t^{ \lambda + \frac{1+j}{k+2}},
$$ o\`u  $\tilde \delta_{i,u}(t)$ est l'image du cycle $ \delta_{i,u}(t)$
par le changement du variable $z \rightarrow Z_u.$

Il faut remarquer  que les coefficients du changement du
variable $w_{i,u}$ sont  des fonctions quasihomog\`enes  de  
$q_{0,u},$  $  \cdots,$  $
q_{\mu -k -j,u}.$ 
Notamment
$$ w_{i,u}(\tau^{\mu -k-1}q_{0,u},\tau^{\mu -k-2}q_{1,u}, \cdots,
\tau q_{\mu -k-2 ,u})= \tau^{- i } w_{i,u}(q_{0,u}, \cdots,
q_{\mu -k-2 ,u}).$$
On peut d\'eduire de la quasihomog\'en\'eit\'e,
$$ w_{i,u}(q) = \omega^{-u i}w_{i,0}(q).$$
Par cons\'equence, les coefficients du d\'eveloppement de $I^{(i)}_u(t)$
aussi satisfont la relation analogue:
$$ v^{(i)}_{j,u} = (\omega^u)^{\beta i}v^{(i)}_{j,0},$$
pour certain $\beta \in {\bf Q}.$
Ce dernier entra\^{i}ne que
$$ (\rho_u^{(1)}, \cdots, \rho_u^{(k+1)}) = (\rho_0^{(1)}, \cdots,
\rho_0^{(k+1)})$$ pour tous les $u, 1 \leq u \leq k+1.$

Ceci ach\`eve  la d\'emonstration de l'\'etape 2.

{\bf  Le  cas  ${\mathcal   P}^{(\ell)}_{x^0}(s,\frac{\partial}{\partial
s}),$
$\ell \geq 1.$}
Pour voir que le terme d'ordre maximal $\mu +1 $ de l'op\'erateur
${\mathcal  P}^{(\ell)}_{x^0}(s,\frac{\partial}{\partial s})$ est \'egal
\`a $ (s_0-{\tilde s}_0)\Delta_\mu(s),$ il  suffit  de  remarquer
dans la Proposition ~\ref{prop3}, ii),
$$ det {\tilde \Sigma} = (F(x^0,s')+ s_0) \Delta_\mu(s) =
(s_0-{\tilde s}_0)\Delta_\mu(s).$$
On peut voir d'apr\`es la construction de l'op\'erateur ${\mathcal
P}^{(\ell)}_{x^0}(s,\frac{\partial}{\partial s})$  que  la longueur de 
l'indice  $\mid  \alpha
\mid$ prend  uniquement les deux valeurs 0 et 1. La partie
$ det (\tilde S \frac {\partial}{\partial s_0} - \tilde L - \tilde 
V(s'))$  produit  les
op\'erateurs avec $ \alpha =0,$ pourtant
les  $\tilde T'_r (s,\frac{\partial}{\partial s})$ produisent ceux avec
$\alpha = {\bf e}_i = (0,\cdots, 0, \check{1}^{ i} ,0, \cdots,0).$
Une \'egalit\'e evidente se d\'eduit de la m\^eme Proposition,
$$ \{\mbox{ mineurs}\; (j \times j)\;\mbox{ de}\; S'(s) \}  =  \{
\mbox{ mineurs}\; (j \times j)\;\mbox{ de}
S(s) \}, 1 \leq j \leq \mu+1.$$
Donc, le raisonnement sur la multiplicit\'e du facteur
$(s_0 -  s_{0,1}(s'))^{k-j  +1}$  dans  les  coefficients  
$\tilde P_{\mu-j, 0, x^0}^{(\ell)}(s),$ 
est parall\`el \`a celui de $P_{\mu-j}.$ 

Quant \`a la multiplicit\'e de
$ \tilde P_{\mu -j,{\bf e}_i, x^0}^{(\ell)}(s), $  elle est d\'etermin\'ee
par les coefficients de $\tilde T'_r (s, \frac  {\partial}{\partial  s})$
i.e. ceux de  l'op\'erateur   $\tilde T_r (s,  \frac  {\partial}{\partial
s_0})$ introduit dans (2.9)'. De cette d\'efinition
il est possible de voir 
que
$$ \tilde P_{\mu -j,{\bf e}_i,x^0}^{(\ell)}(s) = 
(s_0  -  s_{0,1}(s'))^{k-j  +1}(\mbox{polyn\^ome
de degr\'e } \mu -k -1  \;\mbox {en } s_0), 1 \leq j \leq k+1, $$
sur la strate $D_M^{(k)}.$

C'est  \`a  dire,  un  point  de
$(s', s_0(s') ) \in D_M^{(k)}$ appartient  \`a la   strate $\tilde D_M^{(k)}$ 
 $= D_M^{(k)}
\cap \{s_0 =\tilde {s}_0 \} $ si et seulement si
$(s_0 - {\tilde s}_0  )^{k-j+2}  \mid  P_{\mu -j +1, \alpha, x^0
}^{(\ell)}(s',s_0)$ pour $s'$ t.q. $s_0(s') ={\tilde s}_0,$ et
$1\leq j \leq k+1.$ Le r\'esultat sur la strate $D_{M,-}^{(k)}$
s'en d\'eduit aussi imm\'ediatement.

  Quant \`a l'\'enonc\'e $ii),$  il suffit  d'\'ecrire  de  nouveau
l'\'equation d\'eterminante comme $(3.6),$ $ (3.7), $ et appliquer
le m\^eme argument. C.Q.F.D.

\begin{remark}

{\em i) Soit $s_0=$ ${\tilde s}_0 = -F(x^0,s')$  un  point  singulier  de
l'op\'erateur ${\mathcal  P}^{(\ell)}_{x^0}(s,\frac{\partial}{\partial s})$
t.q.  $(s',  {\tilde  s}_0)  \in  D_M^{(-1)}.$  En  d\'epit  de  la
pr\'esence  de  singularit\'es,  pour  $\Delta_j  =   \xi^{(j)}   -
\xi^{(j-1)}, j = 2, \cdots \mu +1, $  un
cycle  \'evanescent  de  $F(x,s')  +s_0  =0,  $  $K^\lambda_{k,  x^0,
\Delta_j}(s) $ est holomorphe au point $ s_0 = {\tilde  s}_0.$
Autrement dit, $K^\lambda_{k, x^0, \gamma^j}(s) $ ne  se  ramifie
qu'aux $\mu$ points de  $D_M.$  Cela  d\'ecoule  du  fait  que
l'indice d'intersection entre le cycle $\gamma^0$ introduit  dans
Remarque ~\ref{remark2} iii) et $\gamma^j$ est \'egal \`a $0.$
}

{\em ii) Si on obtient $ \rho_1, \cdots, \rho_{\mu+1}$ comme
exposants
caract\'eristiques de ${\mathcal P}^{(\ell)}_{x^0}(s,\frac
{\partial}{\partial s})$  \`a $s_0 = s_0(s'),$
il existe  un cycle  $\gamma_i \in H_1(X^{(\mu,\nu)}_s, {\bf C}),$ tel que
$$ \int_{ Reg(\gamma_i)} (x -x^0)^k (F(x,s') + s_0)^\lambda dx \sim
\sum_{j \geq 0}
a^{(k)}_{i,j}(s_0 -s_0(s') )^{\rho_i +j}.
$$ }
\label{remark5}
\end{remark}


{\center{\section{\bf
Calcul des exposants caract\'eristiques
}}}

D\`es  qu'on  a  \'etablit  le   th\'eor\`eme   d'isomonodromie
(Th\'eor\`eme ~\ref{thm3}) le  long des strates
$D_{M,-}^{(k)},{\tilde D}_{M}^{(k)}  $ il est possible de
calculer  l'\'equation
d\'eterminante  des  exposants  pour   $s'   \in   D_{M,-}^{(k)}$  etc.
 \`a
l'autre point ${s'}^{\bullet} \in D_{M,-}^{(k)}$ etc.  pour lequel  le
calcul  se facilite.
 Le lieu de ramification de l'int\'egrale $K_0^{\lambda}(s)$
est simplement stratifi\'e par les strates $  D_{M,-}^{(k)}.  $
Quant \`a $K^{\lambda}_{k, x^0, \gamma}(s),$ son lieu de ramification
a des strates
$D_{M,-}^{(k)}, {\tilde D}_{M}^{(k)}$.
D\`es que on a pos\'e la condition que $F(x,s')$ de (1.1) n'ait que
la singularit\'e du type de Morse, on se restreint aux strates de types
$D_{M,-}^{(0)},{\tilde D}_{M}^{(0)}$ qui peuvent \^etre represent\'ees
par les points suivants:
$$ (s_0, s_1, 0, \cdots,0) \in D_{M,-}^{(0)}, \;\;, s_1 \not = 0, \Delta_\mu
(s_0, s_1, 0, \cdots,0) =0,$$
$$ (s_0,0, s_2, 0, \cdots,0) \in {\tilde D}_{M}^{(0)}, \;\;, s_2 \not =
0, \Delta_\mu (s_0,0, s_2, 0, \cdots,0) =0.$$
Ici on remarque que $(s_0,0, s_2, 0, \cdots,0)$ appartient \`a  l'ensemble 
de bifurcation si $\mu$ est impair. Ce cas sera donc exclu dans les calculs
suivants.

La propostion ci-dessus sert de base \`a l'analyse asymptotique de
l'int\'egrale
$K_{k, x^0}^{\lambda}(s)$ pr\`es de ses points de ramification.
Par   la   suite    on    notera    $K_{k}^{\lambda}(s)    =K_{k,
0}^{\lambda}(s),$ avec $x^0 =0.$
\begin{prop}
L'int\'egrale de p\'eriodes $K_0^{\lambda}(s_0, s_1, 0, \cdots,0)$
(not\'e desormais $K_0^{\lambda}(s_0,  s_1)$)  satisfait  
l'\'equation du type de Fuchs de degr\'e $\mu.$
$$
[\prod _{j = 0}^{\mu -1 }
(\vartheta_{s_0}-{\ell}_j  + j -\lambda)
+ \frac{s_1}{\mu +1} \psi^{\mu}] K_0^{\lambda}(s_0, s_1)=0
\leqno(4.1).$$
En  revanche,    les   int\'egrales   $K_k^{\lambda}(s_0,   s_1),$
d\'efinies   parall\`element   \`a   $K_0^{\lambda}(s_0,   s_1),$
satisfont les \'equations suivantes:
$$
[(\vartheta_{s_0}-{\ell}_{k-1}  + \mu -1 - \lambda)
\prod _{j = 0}^{\mu -1 }
(\vartheta_{s_0}-{\ell}_{k+j}  + j -\lambda)
+ \frac{s_1}{\mu +1} \psi^{\mu}
(\vartheta_{s_0}-k -1 - \lambda)
 ] K_k^{\lambda}(s_0, s_1)=0.
\leqno(4.2)$$

Supposons  $\mu = \ 2m+2$.
Alors les op\'erateurs diff\'erentiels qui annulent
$K^{\lambda}_{j} (s_0,0, s_2, 0)  $, $0 \leq j \leq\mu -1  $,
sont d'ordre  $(\mu +1):$ 

$
S_{0,j-1}(\vartheta_{s_0} + \mu -1- j,\lambda) T_{0, m-1}
(\vartheta_{s_0} + m -j +1, \lambda) \, S_{j,m}(\vartheta_{s_0} -j, \lambda)
(\vartheta_{s_0} - (\lambda +\frac{\mu}{\mu +1})
\\
\qquad\qquad\qquad\qquad\qquad\qquad\qquad\qquad\qquad
+ \mu -1 -j) K_{2j}(s_0,0, s_2, 0)
$
$$
= (\frac{2s_2}{\nu+2})^2 \phi^{2m+1} ( \vartheta_{s_0} -
(\lambda +\frac12+j))
K_{2j}(s_0,0, s_2, 0) , \leqno (4.3)_{2j,o}
$$

$ T_{0,j-1}( \vartheta_{s_0} + \nu - j,\lambda)
S_{0,m-1}(\vartheta_{s_0}    +    m-j    +1,    \lambda)T_{     j
,m}(\vartheta_{s_0} -m, \lambda)
(\vartheta_{s_0} -(\lambda + \frac{\mu -1}{\mu +1})
\\
\qquad\qquad\qquad\qquad\qquad\qquad
  + \mu  -j -2)K_{2j+1}(s_0,0, s_2, 0)  $
$$
= (\frac{2s_2 }{\mu +1}) \phi^{2m+1}(\vartheta_{s_0} - (\lambda + j  +m
+ \frac{3}{2}))
K_{2j+1}(s_0,0, s_2, 0) .\leqno(4.3)_{2j+1,o}
$$
\par\noindent
Ici on a adopt\'e la notation
$$
\vartheta_{s_0} =  s_0 \frac{\partial}{\partial s_0}, \;\;\;
\psi = - \frac{\mu }{\mu +1} s_1 \frac{\partial}{\partial s_0},
\;\;
\phi = - \frac{\mu-1 }{\mu +1} s_2 \frac{\partial}{\partial s_0},
\;\; \ell_j = \frac{j+1}{\mu +1} .$$

$
S_{\alpha,\beta}(X,\lambda) =
\prod _{\ell = \alpha}^{\beta}
(X - (\lambda+ \frac{-(\mu -1) \ell + 1}{ \mu +1}) ),\\
 T_{\alpha,\beta}(X,\lambda) = \prod _{\ell = \alpha}^{\beta}
(X - (\lambda+ \frac{2 -(\mu -1) \ell }{\mu +1}) ).\\$
\label{prop8}
\end{prop}

{\bf D\'emonstration}
Le principe du calcul s'appuie soit sur la Proposition 2.4,  soit
sur un calcul pareil \`a celui de  \cite{AT1},  Prop.  2.2,  2.4.
Voir aussi \cite{AT2}.

Il est facile de voir qu'un changement  des variables,
$$ t_0 = \frac{s_0}{\mu} , \;\;\; t_1  =  -  \frac{s_1}{\mu  +1},
 $$
transforme les \'equations dans une \'ecriture plus  l\'eg\`ere
des \'equations $(4.1), (4.2).$
$$
[\prod _{j = 0}^{\mu -1 }
(\vartheta_{t_0}-{\ell}_j  + j -\lambda)
 - t_1^{\mu +1} (\frac{\partial}{\partial t_0}  )^{\mu}]
K_0^{\lambda}(t_0, t_1)=0
\leqno(4.1)'.$$

$$
[(\vartheta_{t_0}-{\ell}_{k-1}  + \mu -1 - \lambda)
\prod _{j = 0}^{\mu -1 }
(\vartheta_{t_0}-{\ell}_{k+j}  + j -\lambda)
 - \frac{t_1^{\mu +1}}{\mu +1}(\frac{\partial}{\partial t_0}) ^{\mu}
(\vartheta_{t_0} -k -1 - \lambda)
 ] K_k^{\lambda}(t_0, t_1)=0
\leqno(4.2)'.$$

Ici on  a  not\'e  $K_k^{\lambda}(t_0,  t_1),$ les int\'egrales de
p\'eriodes,  apr\`es  le  changement  de
variables ci-dessus  (un abus de notation).

En suite,
on  calcule  les  exposants  caract\'eristiques  des  \'equations
$(4.1)', (4.2)'$ \`a leurs points singuliers.
D'abord on d\'etermine les points singuliers de ces   \'equations.
Dans ce but regardons le symbole  principal des op\'erateurs.

On remarque facilement
$$ \mbox{  le  symbole  principal  de  (4.1)' est}  (t_0^{\mu  }
-1)\sigma ( \frac{\partial}{\partial t_0} )^{\mu}$$
$$ \mbox{  le  symbole  principal  de  (4.2)' est } t_0 (t_0^{\mu }
-1)\sigma ( \frac{\partial}{\partial t_0} )^{\mu +1}$$
si on fixe $t_1$ \`a une valeur non nulle convenablement choisie.
Donc les points singuliers r\'eguliers de (4.1)'  =  $\{  1,
\omega, \omega^2, \cdots, \omega^{\mu -1}, \infty \},$ celles  de
(4.2)' =   $\{  0,  1,  \omega,  \omega^2,  \cdots,  \omega^{\mu
-1},\infty \}$ avec $\omega = e ^{ \frac{2\pi i}{\mu}}.$
Gr\^ace   \`a la    remarque   ~\ref{remark5}    i),    l'int\'egrale
$K^{\lambda}_{k,\gamma } (t_0,t_1)$  se  ramifie  hors  de   $\{t_0
=0\}$ pour   $\gamma$ le cycle \'evanescent en question.
\begin{prop}

Les exposants caract\'eristiques de (4.1)' au point singulier
$t_0 = \omega^j, 0 \leq j \leq \mu -1, $ (correspondant \`a un point de la
strate $ D^{(0)}_{M,-}$) sont comme suit:
$$ \rho = 0, 1, \cdots, \mu -2, \frac {1}{2} + \lambda. $$

Les exposants caract\'eristiques de (4.2)' au point singulier
$t_0 = \omega^j, 0 \leq j \leq \mu -1,$ (correspondant \`a $ D^{(0)}_{M,-}$):
$$ \rho = 0, 1, \cdots, \mu -1, \frac {1}{2}  + \lambda. $$
Au point $t_0 =0,$
$$ \rho = 0, 1, \cdots, \mu -1, \lambda + k+1. $$
Au point singulier
$t_0 = +\infty,$:
$$ \rho =\frac {\mu j - (k+1)}{\mu+1}- \lambda, \;\;\; 0 \leq j
\leq \mu. $$
Les exposants sont rang\'es de sorte que leur somme soit  \'egale
\`a  $\frac{\mu (\mu +1)^2}{2}$ (relation de Riemann-Fuchs).

Les exposants caract\'eristiques de $(4.3)_{2j,o}$ au point singulier
$t_0 = 0,$(correspondant \`a $ {\tilde D}^{(0)}_{M}$) :
$$ \rho = 0, 1, \cdots,2m, \lambda + j +\frac {1}{2} . $$

Les exposants caract\'eristiques de $(4.3)_{2j+1,o}$ au point singulier
$t_0 = 0,$(correspondant \`a ${\tilde D}^{(0)}_{M}$ ) :
$$ \rho = 0, 1, \cdots,2m, \lambda + j + m  + \frac{3}{2}. $$

\label{prop9}
\end{prop}

{\bf D\'emonstration}

La preuve s'appuie sur les  calculs formels d'apr\`es les \'equations
$(4.1)',$ $ (4.2)',$ $ (4.3)$
et la d\'efinition des exposants caract\'eristiques.

Pour obtenir l'\'equation relative de $(4.2)'$ aux points $t=
t_0 = \omega^j,$ il faut savoir compter les
exposants caract\'eristiques de l'op\'erateur
$$ (\vartheta_t + \alpha_0) \cdots (\vartheta_t + \alpha_{\mu})
- \partial_t^{\mu}(\vartheta_t + \gamma)$$
$$ = t(t^{\mu} -1) \partial_t^{\mu +1} + [(\frac{\mu(\mu +1)}{2} + \alpha_0
+ \cdots + \alpha_{ \mu})t^{\mu} -\gamma- \mu  ]\partial_t^{\mu} + \mbox{
des termes de bas ordre. }$$
La seconde \'egalit\'e d\'ecoule des formules:

$$ (\vartheta_t + \alpha_0) \cdots (\vartheta_t + \alpha_{\mu})=
t^{\mu +1}  \partial_t^{\mu +1} +  (\frac{\mu(\mu +1)}{2} + \alpha_0
+ \cdots + \alpha_{ \mu})t^{\mu}\partial_t^{\mu} + {\cdots}$$
$$ \partial_t^{\mu}(\vartheta_t + \gamma)= (\vartheta_t + \gamma + \mu)
\partial_t^{\mu}.$$
En appliquant la d\'efinition ~\ref{dfn31}, on obtient l'\'equation relative
$$ \rho (\rho -1)\cdots (\rho - \mu +1)
(\mu (\rho - \mu) + \frac{\mu(\mu +1)}{2} +
\sum_{ j=0}^{\mu } \alpha_j -\gamma -\mu )=0.$$
Afin de s'adapter \`a $(4.2)',$  on met $ \alpha_j = \frac{\mu}{\mu+1}j
-\lambda - \frac {k+1}{\mu +1},\;$ $0 \leq j \leq \mu -1, $
$\alpha_{\mu}= -\lambda +
\mu -1 -\frac {k}{\mu +1},$ $ \gamma = - \lambda - k -1.$ En somme on obtient
les
exposants caract\'eristiques \'enonc\'es.

D'une fa\c{c}on analogue, on calcule les exposants caract\'eristiques
de $(4.1)'$ \`a
$ t = \omega^j$ \`a l'aide de
$$ (\vartheta_t + \alpha_0) \cdots (\vartheta_t + \alpha_{\mu-1})
- \partial_t^{\mu}$$
$$ = (t^{\mu} -1) \partial_t^{\mu} + (\frac{\mu(\mu -1)}{2} + \alpha_0
+ \cdots + \alpha_{ \mu-1})t^{\mu -1}\partial_t^{\mu -1} + \mbox{
des termes de bas ordre}.$$
Dans ce cas l'\'equation relative est
$$ \rho (\rho -1)\cdots (\rho - \mu + 2)
(\mu (\rho - \mu +1 ) + \frac{\mu(\mu - 1)}{2}+\sum_{ j=0}^{\mu-1 }
\alpha_j)=0.$$ Tous les autres cas se d\'eduisent imm\'ediatement
de $(4.3).$
C.Q.F.D.

L'\'enonc\'e suivant d\'ecoule d'un r\'esultat classique de 
Picard \cite{Pic}, \cite{Ph}.

\begin{prop}
i) Soit $\Delta_j$ le cycle \'evanescent qui correspond \`a la valeur
critique $ t = \omega^j$ de l'int\'egrale $K_0^\lambda(t)$. Soit 
$\Delta_{j+1}$ un autre cycle tel que $\langle \Delta_j,\Delta_{j+1} \rangle 
\not= 0.$ Alors la matrice de monodromie de Picard-Lefschetz autour du point singulier $ t_0 = \omega^j,$ s'\'ecrit sous la forme suivante:
$$\left[
\begin{array}{c}
 K_{k, \Delta_j}^\lambda ( e^{2\pi i }(t-t_0) + t_0)\\
 K_{k, \Delta_{j+1}}^\lambda ( e^{2\pi i }(t-t_0) + t_0)\\
\end{array}
\right]
=
\left[
\begin{array}{cc}
e^{2\pi(\lambda + \frac{1}{2})i} & 0 \\
 e^{2\pi(\lambda + \frac{1}{2})i} & 1\\
\end{array}
\right] \left[
\begin{array}{c}
 K_{k, \Delta_j}^\lambda ( t)\\
 K_{k, \Delta_{j+1}}^\lambda ( t)\\
\end{array}
\right],
 \leqno(4.4)$$ 
Pour les int\'egrales de $(4.2)', (4.3)_{2j,o}, (4.3)_{2j+1,o}$
la matrice de monodromie autour du point $t=0$ s'\'ecrit comme
(4.4).

ii) On a comportement asymptotique pr\`es de $t=\omega^j$ comme suit.

Si $\lambda + \frac{1}{2} \notin \bf Z$
$$ K_{k, \Delta_j}^\lambda (t) = A_k(t-t_0)^{\lambda + \frac{1}{2}}(1+
hol(t-t_0)), \leqno(4.5)$$
$$ K_{k, \Delta_{j+1}}^\lambda (t) = A_k(t-t_0)^{\lambda + \frac{1}{2}}(1+
hol(t-t_0)) + hol(t-t_0),$$
avec $A_k \not =0.$

Si $\lambda + \frac{1}{2} \in \bf Z,$ on a le comportement asymptotique 
(4.5) pour $ K_{k, \Delta_j}^\lambda (t),$ et
$$ K_{k, \Delta_{j+1}}^\lambda (t) = A_k(t-t_0)^{\lambda + \frac{1}{2}} \log
(t-t_0)(1+hol(t-t_0)) + hol(t-t_0). \leqno(4.6)$$

iii) Pour $\mu$ pair, on a comportement asymptotique pr\`es de $t=0$ 
comme suit.

Si $k$ pair et  $\lambda + \frac{1}{2} \notin \bf Z,$
$$ K_{k, \Delta_j}^\lambda (t) = A_k t^{\lambda + \frac{1+k}{2}}(1+
hol(t)), \leqno(4.7)$$
$$ K_{k, \Delta_{j+1}}^\lambda (t) = A_k t^{\lambda + \frac{1+k}{2}}(1+
hol(t)) + hol(t).$$

Si $k$ impair  et  $\lambda + \frac{1}{2} \notin \bf Z,$
$$ K_{k, \Delta_j}^\lambda (t) = A_k t^{\lambda + \frac{\mu + k}{2}}(1+
hol(t)), \leqno(4.8)$$
$$ K_{k, \Delta_{j+1}}^\lambda (t) = A_k t^{\lambda + \frac{\mu +k}{2}}(1+
hol(t)) + hol(t).$$
Dans les cas o\`u  $\lambda + \frac{1}{2} \in \bf Z,$ il faut modifier le d\'eveloppement
asymptotique de $K_{k, \Delta_{j+1}}^\lambda (t)$ comme dans
(4.6).
\label{prop43}
\end{prop}


 \vspace{2pc}
{\center{\section{\bf
D\'emonstration du Th\'eor\`eme ~\ref{thm0} }}}

Ici on d\'emontre le th\'eor\`eme suivant qui implique 
immediatement le th\'eor\`eme~\ref{thm0}, si on met $\nu =2.$

\begin{thm} Dans la situation ci-dessus, on consid\`ere
l'int\'egrale
hyperelliptique 
$I_{P_{K, m}}(s)$
prise le long d'un cycle $\gamma_s =
\{(x,y)\in {\bf R}^2: H(x,y)+ s_0=0 \},$
$$   I_{P_{K,m}}(s)=   \int_{\gamma_s}P_{K,m}(x,y)    dx    ,
\leqno(5.1)$$
avec $ K \in {\bf N}, m \in {\bf Z}.$ Supposons $ I_{P_{K, m}}(s)\not \equiv
0.$  Alors on a r\'esultats suivants.

i) Si $\mu$ pair, la multiplicit\'e $N$ des z\'eros
de  
l'int\'egrale $I_{P_{K, m}}(s)$ \`a l'un des points de ramification
${\tilde t}\in \{ s_0^{(1)}, \cdots ,s_0^{(\mu)}    \}$ v\'erifie:
$$ N \leq 
2[\frac{m}{\nu} + \frac{K+\mu}{2}]. \leqno(5.2)$$

ii) Supposons que $k_1=0$ dans l'expression (1.3).
Alors la multiplicit\'e $N$ des z\'eros
de  l'int\'egrale $I_{P_{K, m}}(s)$ avec $\mu$ impair, 
\`a l'un des points de ramification
${\tilde s_0}\in \{ s_0^{(1)}, \cdots ,s_0^{(\mu)}    \}$ 
v\'erifie:
$$ N \leq max \{ \mu-1, 
2[\frac{m}{\nu} + \frac{3}{2}] \}.
\leqno(5.3)$$

iii) La multiplicit\'e $N$ des z\'eros
de  
l'int\'egrale $I_{P_{K, m}}(t)$ au point 
$\tilde s_0 \not \in \{ s_0^{(1)}, \cdots ,s_0^{(\mu)}\}$ ne 
d\'epasse pas $\mu +K .$
\label{thm51}
\end{thm}

Pour  d\'emontrer  le Th\'eor\`eme  ~\ref{thm51},  il   faut   d\'ej\`a
quelques lemmes \'el\'ementaires. Les notations sont
celles des chapitres pr\'ec\'edents. 

\begin{lem}
Soient   $x^{(2)},   \cdots,   x^{(L)}$   les   points
non-critiques    de
l'application $x \rightarrow F(x,s').$ Alors l'\'egalit\'e suivante a
lieu entre les deux int\'egrales
$$I_{P_{K,m}(x,y),\Gamma}(s) = ( (\frac{\partial}{\partial s_0})^{-1}
\frac{\partial}{\partial f_1(x^{(2)},s')  })^{k_2}\cdots
( (\frac{\partial}{\partial s_0})^{-1}
\frac{\partial}{\partial  f_1(x^{(L)},s')
})^{k_L}I_{(x - x^{(1)})^{k_1}y^m,\Gamma}(s).$$
Ici
$$\frac{\partial}{\partial f_1(x^{(i)},s')  }= \frac{1}{\mid
grad _{s'} f_1(x^{(i)},s') \mid } \sum_{j=1}^{\mu -1}
\frac{\partial}{\partial s_j} f_1(x^{(i)},s')
\frac{\partial}{\partial s_j},
$$ et
$$ f_1(x^{(i)},s')  = \frac{\partial}{\partial x} F(x,  s')\mid_{x  =
x^{(i)}}.$$
O\`u l'action de $(\frac{\partial}{\partial s_0})^{-1} $
sur la fonction de 
Dulac   d\'efinie par l'expansion asymptotique
pr\`es du point $s_0 = x$ dans un secteur simplement connexe du plan
complexe $ {\bf C}_{s_0}$ est donn\'ee d'une fa\c{c}on suivante:
$$ (\frac{\partial}{\partial s_0})^{-1} [(s_0- x)^{\rho}
(\sum_{j \geq 0}a_j (s_0- x)^{j}) ] =
(s_0- x)^{\rho}
(\sum_{j \geq 0}\frac {a_j}{ \rho + j +1} (s_0- x)^{j +1 }) .$$

\label{lem51}
\end{lem}
{\bf D\'emonstration }

On a le d\'eveloppement analogue \`a $(2.3)'$  en 
$(x - x^{(i)}),$
$$
F(x,s') = (x - x^{(i)})^{\mu + 1}+ \sum_{\ell =1}^{\mu }
f_{\ell}(x^{(i)}, s') (x - x^{(i)})^{ \ell } + F (x^{(i)}, s').$$

De celui-ci, il est facile de voir
$$\frac{\partial}{\partial f_1(x^{(i)},s')  } \int_{\Gamma}
(x - x^{(1)})^{k_1} (F(x,s') + s_0)^{\lambda} dx
=\frac{\partial}{\partial s_0}
\int_{\Gamma}(x - x^{(i)})  (x - x^{(1)})^{k_1} (F(x,s') +s_0
)^{\lambda} dx , \;\;\;\; i  = 2,3,\cdots.
$$
De ce dernier, on d\'eduit la relation d\'esir\'ee. C.Q.F.D.

\begin{lem}
Si l'int\'egrale $I_{(x - x^{(1)})^{k_1}y^m,\Gamma}(s)$ admet  un
d\'eveloppement asymptotique pr\`es  du  point  de  ramification  $s_0  =
s_0(s'),$ $(s',s_0(s')) \in D_M^{(0)},$
$$I_{(x - x^{(1)})^{k_1}y^m,\Gamma}(s)\sim
\sum_{j \geq 0}a_j^{(k_1)}(s') (s_0- s_0(s') )^{\rho + j}+
\sum_{j \geq 0}a_{j,1}^{(k_1)}(s')(s_0- s_0(s'))^{\rho+j}\log (s_0- s_0(s')) 
,$$
avec $a_0^{(k_1)}(s')\cdot a_{0,1}^{(k_1)}(s') \not = 0$ sur la strate
$D_{M,-}^{(0)}$ (ou $D_{M,1}^{(0)}$ ), alors d\'eveloppement
asymptotique ci-dessous a lieu:
$$I_{(x - x^{(1)})^{k_1} (x -x^{(i)})y^m,\Gamma}(s)
\sim \sum_{j \geq 0}a_{i,j}^{(k_1)}(s') (s_0- s_0(s') )^{\rho + j}+
\sum_{j \geq 0}a_{i,j,1}^{(k_1)}(s')(s_0- s_0(s'))^{\rho+j}\log (s_0- s_0(s')) 
,$$
avec $ a_{i,0}^{(k_1)}(s')\cdot a_{i,0,1}^{(k_1)}(s')\not = 0,$ sur 
la strate $D_{M,-}^{(0)}$
(ou $\tilde D_{M}^{(0)}$).
\label{lem52}
\end{lem}

{\bf D\'emonstration }

D'apr\`es le Lemme~\ref{lem51}, on a d\'eveloppement asymptotique
suivant:
$$
( (\frac{\partial}{\partial s_0})^{-1}
\frac{\partial}{\partial f_1(x^{(i)},s')  })
I_{(x - x^{(1)})^{k_1}y^m,\Gamma}(s) $$
$$\sim a_0^{(k_1)}(s')(-\frac{\partial}{\partial    f_1(x^{(i)},s')
}s_0(s')) (s_0- s_0(s') )^{\rho} +
a_{0,1}^{(k_1)}(s')(-\frac{\partial s_0(s')}{\partial f_1(x^{(i)},s')})
(s_0- s_0(s'))^{\rho}\log (s_0- s_0(s')) +$$
$$+\frac{1}{\rho +1}(\frac{\partial}{\partial f_1(x^{(i)},s')
}a_0^{(k_1)}(s'))(s_0- s_0(s'))^{\rho +1} + \cdots.$$
Donc l'\'enonc\'e d\'esir\'e est r\'eduit \`a l'in\'egalit\'e,
$$\frac{\partial}{\partial    f_1(x^{(i)},s') }s_0(s') \not =0.$$

  Il est facile de voir que $\frac{\partial}{\partial
f_1(x^{(i)},s') },$ consid\'er\'e  comme  un  champ  de  vecteur  dans
${\bf  C}^{\mu}_s,$  ne  peut  pas  \^etre  un   champ   de   vecteurs
logarithmique  i.e.  $\frac{\partial}{\partial     f_1(x^{(i)},s')  }
\not \in Der_{{\bf C}^{\mu}_s}(log\;D).$ Cela  d\'ecoule 
du  fait
que les \'el\'ements  de  $Der_{{\bf  C}^{\mu}_s}(log\;D)$  sont  \`a
coefficients        polynomiaux        quasihomog\`enes         (voir
la Proposition~\ref{prop2},   le  lemme~\ref{lem1})     avec     des     poids
quasihomog\`enes positifs, par contre
$\frac{\partial}{\partial      f_1(x^{(i)},s')   }$    est    \`a
coefficients constants de  poids  negatif.  C'est- \`a -dire que  le
vecteur normal \`a la surface $s_0 = s_0(s')$ (i.e.  $D_M^{(0)}$  )  dans
${\bf C}^\mu_s$ n'est pas normal \`a $\frac{\partial}{\partial
f_1(x^{(i)},s')}.$
On en d\'eduit que
$$ \frac{1} {\mid grad_{s'} f_1(x^{(i)},s')  \mid}
< (0, f_1(x^{(i)},s')  ), (1, grad_{s'}s_0(s'))> =
\frac{\partial}{\partial     f_1(x^{(i)},s')  }s_0(s') \not = 0.$$
C.Q.F.D.

En r\'ep\'etant l'argument du Lemme~\ref{lem52}, pour chaque point de
ramification $ (s',s_0(s')) \in      D_M^{(0)},$ on obtient l'\'enonc\'e 
suivant.
On notera dans la suite la strate $  D_M^{(0)}  \cap  \{(s',s_0);
s_0 = -F(x^{(1)},s')\}   $   par
$D_{M,1}^{(0)}.$ C'est une  strate
analogue    \`a     ${\tilde D}_M^{(0)}.$

\begin{prop}
La multiplicit\'e des z\'eros de $I_{(x - x^{(1)})^{k_1}y^m,\Gamma}$ $ (s) 
\not \equiv 0$ est \'egale  \`a celle des z\'eros de l'int\'egrale
$I_{P_{K,m}(x, y),\Gamma} (s)$ \`a chaque point  de  ramification
de la  strate $D_{M,-}^{(0)}$ ou $D_{M,1}^{(0)},$ si
$x^{(2)}, \cdots, x^{(L)}$  ne  sont  pas  des  points  critiques  de
$F(x,s').$
\label{prop53}
\end{prop}
  Pour  savoir  le  comportement  asymptotique  de
l'int\'egrale $I_{(x - x^{(1)})^{k_1} \cdots (x - x^{(L)})^{k_L}
y^m,\Gamma}(s)$ pr\`es de la   strate  $D_{M}^{(0)},$  il  suffit  de
savoir celui de $I_{(x - x^{(1)})^{k_1} y^m,\Gamma}(s).$

{\bf 5.2. D\'emonstration du Th\'eor\`eme~\ref{thm51}}

Pour estimer la multiplicit\'e   des z\'eros de
$I_{P_{K,m}(x, y),\Gamma}(s),$  d'abord on estime  les  exposants
caract\'eristiques aux points de ramification de la strate $D_{M}^{(0)},$
\`a l'aide de la Proposition~\ref{prop9} et de la Proposition~\ref{prop53}.
Dans le cas o\`u $I_{P_{K,m}(x, y),\Gamma}(s)$ est holomorpe \`a
$s_0  =s_0^{(j)}$  pour un certain  $j,$  notre  argument  d'estimation
ci-dessous est toujours valable. Ici on  consid\`ere  le  cas  du
nombre maximal $(= \mu )$ de points de ramification.

La  contribution  maximale $\rho$ des  exposants  caract\'eristiques  \`a
chaque point est classifi\'ee comme suit.

I) cas $\mu$ pair

Au  point  $s_0  =  s_0^{(1)},$  correspondant   au   point   de la   strate
$D_{M,1}^{(0)},$

si $k_1$ pair $\rho = \frac{m}{\nu} + \frac{k_1 + 1}{2}  .$

si $k_1$ impair $\rho = \frac{m}{\nu} + \frac{k_1 + \mu}{2}.$

A chaque point $s_0 = s_0^{(i)},$ $i = 2, \cdots \mu,$
correspondant   au   point   de la   strate $D_{M,-}^{(0)},$
$\rho =\frac{m}{\nu} +{\frac{1}{2}}  .$

II) cas $\mu$ impair

Dans ce cas l\`a, on ne sait
que des exposants caract\'eristiques au point de la
strate
$D_{M,1}^{(0)},$ dont le maximum est toujours \'egal \`a
$\rho =\frac{m}{\nu} +{\frac{1}{2}}  .$

D'apr\`es la D\'efinition ~\ref{dfn11} et la Proposition ~\ref{prop53}, la multiplicit\'e $N$ des z\'eros
de l'int\'egrale peut \^etre \'estim\'e par $2 \rho$ si $\rho$ $\in \bf Z$
et par $[\rho] +1$ si $\rho$ $\notin \bf Z.$

De la  liste ci-dessus se d\'eduit imm\'ediatement l'\'enonc\'e
i), ii) du Th\'eor\`eme ~\ref{thm51}.  

Quant \`a l'\'enonc\'e iii), on utilise le fait que 
$I_{(x - x^{(1)})^{k} y^m,\Gamma}(s)$ satisfait une \'equation 
diff\'erentielle aux singularit\'es r\'eguli\`eres de degr\'e $\mu +1, $ 
${\mathcal P}^{(k)}_{x^{(1)}}(s, \frac {\partial}{\partial s}) $ (voir la
Proposition ~\ref{prop3}, ii)). Cela veut dire que les exposants
caract\'eristiques de 
${\mathcal P}^{(k)}_{x^{(1)}}(s, \frac {\partial}{\partial s}) $ au point
hors discriminant $D$ doivent \^etre $\{0,1,\cdots, \mu \}$
o\`u la solution $I_{(x - x^{(1)})^{k} y^m,\Gamma}(s)$ est holomorphe. 
 D'autre part, le lemme~\ref{lem51} implique que l'exposant
caract\'eristique de  $I_{P_{K,m}(x,  y),\Gamma}(s)$ est sup\'erieur d'au plus 
$k_2+ \cdots + k_L \leq K $  \`a celui de $I_{(x - x^{(1)})^{k} y^m,\Gamma}(s).$
D'o\`u vient l'\'enonc\'e d\'esir\'e.
{\bf C.Q.F.D.}
                                                                                                      

\vspace{\fill}


%

\noindent

\begin{flushleft}

 \begin{minipage}[t]{6.2cm}

  \begin{center}

{\footnotesize Moscow Independent University\\
Bol'shoj Vlasijevskij Pereulok 11,\\
  MOSCOW, 121002,\\

Russia\\

{\it E-mails}:  tanabe@mccme.ru, tanabesusumu@hotmail.com}

\end{center}

\end{minipage}\hfill

\end{flushleft}

\end{document}